\documentclass[final,3p]{elsarticle}
 \usepackage{graphics}
 \usepackage{graphicx}
 \usepackage{epsfig}
\usepackage{amssymb}
 \usepackage{amsthm}
 \usepackage{lineno}
 \usepackage{amsmath}
   \numberwithin{equation}{section}
\usepackage{mathrsfs}

\journal{Springer} 

\newtheorem{thm}{Theorem}[section]

\newtheorem{lem}[thm]{Lemma}
\newtheorem{prop}[thm]{Proposition}

\biboptions{sort&compress,square}
\allowdisplaybreaks
\begin{document}
\begin{frontmatter}
\author{Tong Wu$^{a}$}
\ead{wut977@nenu.edu.cn}
\author{Yong Wang$^{b,*}$}
\ead{wangy581@nenu.edu.cn}
\cortext[cor]{Corresponding author.}
\author{Sining Wei$^c$}
\ead{weisn835@nenu.edu.cn}
\address{$^a$Department of Mathematics, Northeastern University, Shenyang, 110819, China}
\address{$^b$School of Mathematics and Statistics, Northeast Normal University,
Changchun, 130024, China}
\address{$^c$School of Data Science and Artificial Intelligence, Dongbei University of Finance and Economics, Dalian 116025, P.R.China}

\title{The general Kastler-Kalau-Walze type theorem and the Dabrowski-Sitarz-Zalecki type theorem for odd dimensional manifold with boundary}
\begin{abstract}
In this paper, we give the proof of the general Kastler-Kalau-Walze type theorem and the Dabrowski-Sitarz-Zalecki type theorem on odd dimensional compact manifolds with boundary.
\end{abstract}
\begin{keyword}The Dirac operator; the Kastler-Kalau-Walze type theorem; the Dabrowski-Sitarz-Zalecki type theorem\\

\end{keyword}
\end{frontmatter}
\textit{2010 Mathematics Subject Classification:}
53C40; 53C42.
\section{Introduction}
 Until now, many geometers have studied noncommutative residues. In \cite{Gu,Wo}, authors found noncommutative residues are of great importance to the study of noncommutative geometry. In \cite{Co1}, Connes used the noncommutative residue to derive a conformal 4-dimensional Polyakov action analogy. Connes showed us that the noncommutative residue on a compact manifold $M$ coincided with the Dixmier's trace on pseudodifferential operators of order $-{\rm {dim}}M$ in \cite{Co2}.
And Connes claimed the noncommutative residue of the square of the inverse of the Dirac operator was proportioned to the Einstein-Hilbert action.  Kastler \cite{Ka} gave a
brute-force proof of this theorem. Kalau and Walze proved this theorem in the normal coordinates system simultaneously in \cite{KW} .
Ackermann proved that
the Wodzicki residue  of the square of the inverse of the Dirac operator ${\rm  Wres}(D^{-2})$ in turn is essentially the second coefficient
of the heat kernel expansion of $D^{2}$ in \cite{Ac}.

On the other hand, Wang generalized the Connes' results to the case of manifolds with boundary in \cite{Wa1,Wa2},
and proved the Kastler-Kalau-Walze type theorem for the Dirac operator and the signature operator on lower-dimensional manifolds
with boundary \cite{Wa3}. In \cite{Wa3,Wa4}, Wang computed $\widetilde{{\rm Wres}}[\pi^+D^{-1}\circ\pi^+D^{-1}]$ and $\widetilde{{\rm Wres}}[\pi^+D^{-2}\circ\pi^+D^{-2}]$, where the two operators are symmetric, in these cases the boundary term vanished. But for $\widetilde{{\rm Wres}}[\pi^+D^{-1}\circ\pi^+D^{-3}]$, J. Wang and Y. Wang got a nonvanishing boundary term \cite{Wa5}, and give a theoretical explanation for gravitational action on boundary. In others words, Wang provided a kind of method to study the Kastler-Kalau-Walze type theorem for manifolds with boundary. In \cite{DL}, the authors defined bilinear functionals of vector fields and differential forms, the densities of which yielded the  metric and Einstein tensors on even-dimensional Riemannian manifolds. In \cite{Wu2}, Wu and Wang gave the proof of the Kastler-Kalau-Walze type theorem for the generalized noncommutative residue on 4-dimensional and 6-dimensional compact manifolds with (resp.without) boundary. In \cite{Wu1}, Wu and Wang computed the noncommutative residue $\widetilde{{\rm Wres}}[\pi^+(\nabla_X^{S(TM)}\nabla_Y^{S(TM)}D^{-2})\circ\pi^+(D^{-(n-2)})]$  and $\widetilde{{\rm Wres}}[\pi^+(\nabla_X^{S(TM)}\nabla_Y^{S(TM)}D^{-1})\circ\pi^+(D^{-(n-1)})]$ on even and odd dimensional compact manifolds. Our motivation is to prove a general Kastler-Kalau-Walze type theorem and a general Dabrowski-Sitarz-Zalecki type theorem for odd dimensional manifolds with boundary. That is, we want to compute ${\rm Wres}[\pi^+P_1\circ\pi^+P_2]$, where orders of $P_1,P_2$ are $a_1,a_2$ and $-a_1-a_2+2=n$ for odd dimensional manifolds with boundary. Motivated by \cite{Wu1,Wu2}, we compute the generalized noncommutative residue $\widetilde{{\rm Wres}}[\pi^+(c(X)D^{-1})\circ\pi^+(D^{-(2m-2)})]$, $\widetilde{{\rm Wres}}[\pi^+(\nabla_X^{S(TM)}D^{-1})\circ\pi^+(D^{-(2m-1)})]$ and $\widetilde{{\rm Wres}}[\pi^+(\nabla_X^{S(TM)}D^{-2})\circ\pi^+(D^{-(2m-2)})]$ on odd dimensional manifolds. Our main theorems are as follows.
\begin{thm}\label{thmb111}
Let $M$ be an $n=2m+1$ dimensional oriented
compact spin manifold with boundary $\partial M$, then we get the following equality:
\begin{align}
\label{1111}
&\widetilde{{\rm Wres}}[\pi^+(c(X)D^{-1})\circ\pi^+(D^{-(2m-2)})]\nonumber\\
&=\int_{\partial M} \bigg\{(m-1)\bigg(\partial_{x_n}(X_n)\frac{\pi}{(m+1)!}A_0+X_n\frac{\pi i}{2(m+2)!}A_1+X_nh'(0)\frac{\pi i}{(m+2)!}B_0\bigg) +i(2m^2-m-1)X_nh'(0)\nonumber\\
&\frac{\pi i}{4(m+1)!}C_{0}+i(m^2-2m+1)X_nh'(0)\frac{2\pi i}{(m+2)!}C_1
+i(m-1)X_nh'(0)\frac{\pi i}{2(m+2)!}D_0\bigg\}Vol(S^{n-2})2^m d{\rm Vol_{M}},
\end{align}
where $X=\sum_{j=1}^nX_j\partial_j$ is a vector field, and $A_0$-$D_0$ are defined in section 3.
\end{thm}
\begin{thm}\label{thmb2}
Let $M$ be an $n=2m+1$ dimensional oriented
compact spin manifold with boundary $\partial M$, then we get the following equality:
\begin{align}
\label{1b2773}
&\widetilde{{\rm Wres}}[\pi^+(\nabla_X^{S(TM)}D^{-1})\circ\pi^+(D^{-(2m-1)})]\nonumber\\
&=\int_{\partial M}\bigg\{-\frac{1}{4}Vol(S^{n-2})2^m\partial_{x_n}(X_n)\frac{2\pi i}{(m+2)!}E_0
+\bigg(-\frac{\pi i}{2(m+3)!}E_1+\frac{\pi i}{2(m+3)!}F_0+\frac{\pi }{4(m+3)!}G_0\nonumber\\
&+\frac{\pi }{2(m+2)!}G_1-\frac{(2m^2-m)\pi}{4(m+2)!}H_0
+\frac{((i-1)m^2+mi+3m)\pi i}{(m+3)!}H_1-\frac{(2m^2+3m+1)\pi i}{4(m+2)!}H_2\nonumber\\
&-\frac{(3m+1)\pi i}{(m+2)!}H_3\bigg)Vol(S^{n-2})X_nh'(0)2^m\bigg\}d{\rm Vol_{M}}
\end{align}
where $X=\sum_{j=1}^nX_j\partial_j$ is a vector field, and $E_0$-$H_3$ are defined in section 4.
\end{thm}
\begin{thm}\label{thmb233}
Let $M$ be an $n=2m+1$ dimensional oriented
compact spin manifold with boundary $\partial M$, then we get the following equality:
\begin{align}
\label{3b2773}
&\widetilde{{\rm Wres}}[\pi^+(\nabla_X^{S(TM)}D^{-2})\circ\pi^+(D^{-(2m-2)})]\nonumber\\
&=\int_{\partial M} \bigg\{-(m-1)Vol(S^{n-2})2^m\frac{\partial X_n}{\partial x_n}\frac{\pi i}{(m+1)!}I_0+\bigg((m-1)\frac{\pi i}{2(m+2)!}I_1
+(1-m)\frac{\pi i}{(m+2)!}J_0\nonumber\\
&+\frac{(2m^2+3m-5)\pi }{2(m+1)!}K_0+\frac{(m^2-3m+2)\pi }{(m+2)!}K_1+(m-1)\frac{2\pi i}{(m+2)!}L_0\bigg)Vol(S^{n-2})X_nh'(0)2^m\bigg\}d{\rm Vol_{M}}.
\end{align}
where $X=\sum_{j=1}^nX_j\partial_j$ is a vector field, and $I_0$-$L_0$ are defined in section 5.
\end{thm}
The paper is organized in the following way. In Section \ref{section:2},  we recall some basic facts and formulas about Boutet de
Monvel's calculus and the definition of the noncommutative residue for manifolds with boundary. In Section \ref{section:3},  we prove the general Kastler-Kalau-Walze type theorem for $\widetilde{{\rm Wres}}[\pi^+(c(X)D^{-1})\circ\pi^+(D^{-(2m-2)})]$ on odd dimensional manifolds with boundary. In Section \ref{section:4},  we prove the general Dabrowski-Sitarz-Zalecki type theorem $\widetilde{{\rm Wres}}[\pi^+(\nabla_X^{S(TM)}D^{-1})\circ\pi^+(D^{-(2m-1)})]$ on odd dimensional manifolds with boundary. In Section \ref{section:5},  we prove the general Dabrowski-Sitarz-Zalecki type theorem $\widetilde{{\rm Wres}}[\pi^+(\nabla_X^{S(TM)}D^{-2})\circ\pi^+(D^{-(2m-2)})]$ on odd dimensional manifolds with boundary.
\section{Boutet de Monvel's calculus}
\label{section:2}
 In this section, we recall some basic facts and formulas about Boutet de
Monvel's calculus and the definition of the noncommutative residue for manifolds with boundary which will be used in the following. For more details, see Section 2 in \cite{Wa3}.\\
 \indent Let $M$ be a 4-dimensional compact oriented manifold with boundary $\partial M$.
We assume that the metric $g^{TM}$ on $M$ has the following form near the boundary,
\begin{equation}
\label{b1}
g^{M}=\frac{1}{h(x_{n})}g^{\partial M}+dx _{n}^{2},
\end{equation}
where $g^{\partial M}$ is the metric on $\partial M$ and $h(x_n)\in C^{\infty}([0, 1)):=\{\widehat{h}|_{[0,1)}|\widehat{h}\in C^{\infty}((-\varepsilon,1))\}$ for
some $\varepsilon>0$ and $h(x_n)$ satisfies $h(x_n)>0$, $h(0)=1$ where $x_n$ denotes the normal directional coordinate. Let $U\subset M$ be a collar neighborhood of $\partial M$ which is diffeomorphic with $\partial M\times [0,1)$. By the definition of $h(x_n)\in C^{\infty}([0,1))$
and $h(x_n)>0$, there exists $\widehat{h}\in C^{\infty}((-\varepsilon,1))$ such that $\widehat{h}|_{[0,1)}=h$ and $\widehat{h}>0$ for some
sufficiently small $\varepsilon>0$. Then there exists a metric $g'$ on $\widetilde{M}=M\bigcup_{\partial M}\partial M\times
(-\varepsilon,0]$ which has the form on $U\bigcup_{\partial M}\partial M\times (-\varepsilon,0 ]$
\begin{equation}
\label{b2}
g'=\frac{1}{\widehat{h}(x_{n})}g^{\partial M}+dx _{n}^{2} ,
\end{equation}
such that $g'|_{M}=g$. We fix a metric $g'$ on the $\widetilde{M}$ such that $g'|_{M}=g$.

Let the Fourier transformation $F'$  be
\begin{equation*}
F':L^2({\bf R}_t)\rightarrow L^2({\bf R}_v);~F'(u)(v)=\int_\mathbb{R} e^{-ivt}u(t)dt
\end{equation*}
and let
\begin{equation*}
r^{+}:C^\infty ({\bf R})\rightarrow C^\infty (\widetilde{{\bf R}^+});~ f\rightarrow f|\widetilde{{\bf R}^+};~
\widetilde{{\bf R}^+}=\{x\geq0;x\in {\bf R}\}.
\end{equation*}
\indent We define $H^+=F'(\Phi(\widetilde{{\bf R}^+}));~ H^-_0=F'(\Phi(\widetilde{{\bf R}^-}))$ which satisfies
$H^+\bot H^-_0$, where $\Phi(\widetilde{{\bf R}^+}) =r^+\Phi({\bf R})$, $\Phi(\widetilde{{\bf R}^-}) =r^-\Phi({\bf R})$ and $\Phi({\bf R})$
denotes the Schwartz space. We have the following
 property: $h\in H^+~$ (resp. $H^-_0$) if and only if $h\in C^\infty({\bf R})$ which has an analytic extension to the lower (resp. upper) complex
half-plane $\{{\rm Im}\xi<0\}$ (resp. $\{{\rm Im}\xi>0\})$ such that for all nonnegative integer $l$,
 \begin{equation*}
\frac{d^{l}h}{d\xi^l}(\xi)\sim\sum^{\infty}_{k=1}\frac{d^l}{d\xi^l}(\frac{c_k}{\xi^k}),
\end{equation*}
as $|\xi|\rightarrow +\infty,{\rm Im}\xi\leq0$ (resp. ${\rm Im}\xi\geq0)$ and where $c_k\in\mathbb{C}$ are some constants.\\
 \indent Let $H'$ be the space of all polynomials and $H^-=H^-_0\bigoplus H';~H=H^+\bigoplus H^-.$ Denote by $\pi^+$ (resp. $\pi^-$) the
 projection on $H^+$ (resp. $H^-$). Let $\widetilde H=\{$rational functions having no poles on the real axis$\}$. Then on $\tilde{H}$,
 \begin{equation}
 \label{b3}
\pi^+h(\xi_0)=\frac{1}{2\pi i}\lim_{u\rightarrow 0^{-}}\int_{\Gamma^+}\frac{h(\xi)}{\xi_0+iu-\xi}d\xi,
\end{equation}
where $\Gamma^+$ is a Jordan closed curve
included ${\rm Im}(\xi)>0$ surrounding all the singularities of $h$ in the upper half-plane and
$\xi_0\in {\bf R}$. In our computations, we only compute $\pi^+h$ for $h$ in $\widetilde{H}$. Similarly, define $\pi'$ on $\tilde{H}$,
\begin{equation}
\label{b4}
\pi'h=\frac{1}{2\pi}\int_{\Gamma^+}h(\xi)d\xi.
\end{equation}
So $\pi'(H^-)=0$. For $h\in H\bigcap L^1({\bf R})$, $\pi'h=\frac{1}{2\pi}\int_{{\bf R}}h(v)dv$ and for $h\in H^+\bigcap L^1({\bf R})$, $\pi'h=0$.\\
\indent An operator of order $m\in {\bf Z}$ and type $d$ is a matrix\\
$$\widetilde{A}=\left(\begin{array}{lcr}
  \pi^+P+G  & K  \\
   T  &  \widetilde{S}
\end{array}\right):
\begin{array}{cc}
\   C^{\infty}(M,E_1)\\
 \   \bigoplus\\
 \   C^{\infty}(\partial{M},F_1)
\end{array}
\longrightarrow
\begin{array}{cc}
\   C^{\infty}(M,E_2)\\
\   \bigoplus\\
 \   C^{\infty}(\partial{M},F_2)
\end{array},
$$
where $M$ is a manifold with boundary $\partial M$ and
$E_1,E_2$~ (resp. $F_1,F_2$) are vector bundles over $M~$ (resp. $\partial M
$).~Here,~$P:C^{\infty}_0(\Omega,\overline {E_1})\rightarrow
C^{\infty}(\Omega,\overline {E_2})$ is a classical
pseudodifferential operator of order $m$ on $\Omega$, where
$\Omega$ is a collar neighborhood of $M$ and
$\overline{E_i}|M=E_i~(i=1,2)$. $P$ has an extension:
$~{\cal{E'}}(\Omega,\overline {E_1})\rightarrow
{\cal{D'}}(\Omega,\overline {E_2})$, where
${\cal{E'}}(\Omega,\overline {E_1})~({\cal{D'}}(\Omega,\overline
{E_2}))$ is the dual space of $C^{\infty}(\Omega,\overline
{E_1})~(C^{\infty}_0(\Omega,\overline {E_2}))$. Let
$e^+:C^{\infty}(M,{E_1})\rightarrow{\cal{E'}}(\Omega,\overline
{E_1})$ denote extension by zero from $M$ to $\Omega$ and
$r^+:{\cal{D'}}(\Omega,\overline{E_2})\rightarrow
{\cal{D'}}(\Omega, {E_2})$ denote the restriction from $\Omega$ to
$X$, then define
$$\pi^+P=r^+Pe^+:C^{\infty}(M,{E_1})\rightarrow {\cal{D'}}(\Omega,
{E_2}).$$ In addition, $P$ is supposed to have the
transmission property; this means that, for all $j,k,\alpha$, the
homogeneous component $p_j$ of order $j$ in the asymptotic
expansion of the
symbol $p$ of $P$ in local coordinates near the boundary satisfies:\\
$$\partial^k_{x_n}\partial^\alpha_{\xi'}p_j(x',0,0,+1)=
(-1)^{j-|\alpha|}\partial^k_{x_n}\partial^\alpha_{\xi'}p_j(x',0,0,-1),$$
then $\pi^+P:C^{\infty}(M,{E_1})\rightarrow C^{\infty}(M,{E_2})$. Let $G$,$T$ be respectively the singular Green operator
and the trace operator of order $m$ and type $d$. Let $K$ be a
potential operator and $S$ be a classical pseudodifferential
operator of order $m$ along the boundary. Denote by $B^{m,d}$ the collection of all operators of
order $m$
and type $d$,  and $\mathcal{B}$ is the union over all $m$ and $d$.\\
\indent Recall that $B^{m,d}$ is a Fr\'{e}chet space. The composition
of the above operator matrices yields a continuous map:
$B^{m,d}\times B^{m',d'}\rightarrow B^{m+m',{\rm max}\{
m'+d,d'\}}.$ Write $$\widetilde{A}=\left(\begin{array}{lcr}
 \pi^+P+G  & K \\
 T  &  \widetilde{S}
\end{array}\right)
\in B^{m,d},
 \widetilde{A}'=\left(\begin{array}{lcr}
\pi^+P'+G'  & K'  \\
 T'  &  \widetilde{S}'
\end{array} \right)
\in B^{m',d'}.$$\\
 The composition $\widetilde{A}\widetilde{A}'$ is obtained by
multiplication of the matrices (For more details see \cite{ES}). For
example $\pi^+P\circ G'$ and $G\circ G'$ are singular Green
operators of type $d'$ and
$$\pi^+P\circ\pi^+P'=\pi^+(PP')+L(P,P').$$
Here $PP'$ is the usual
composition of pseudodifferential operators and $L(P,P')$ called
leftover term is a singular Green operator of type $m'+d$. For our case, $P,P'$ are classical pseudo differential operators, in other words $\pi^+P\in \mathcal{B}^{\infty}$ and $\pi^+P'\in \mathcal{B}^{\infty}$ .\\
\indent Let $M$ be a $n$-dimensional compact oriented manifold with boundary $\partial M$.
Denote by $\mathcal{B}$ the Boutet de Monvel's algebra. We recall that the main theorem in \cite{FGLS,Wa3}.
\begin{thm}\label{th:32}{\rm\cite{FGLS}}{\bf(Fedosov-Golse-Leichtnam-Schrohe)}
 Let $M$ and $\partial M$ be connected, ${\rm dim}M=n\geq3$, and let $\widetilde{S}$ (resp. $\widetilde{S}'$) be the unit sphere about $\xi$ (resp. $\xi'$) and $\sigma(\xi)$ (resp. $\sigma(\xi')$) be the corresponding canonical
$n-1$ (resp. $(n-2)$) volume form.
 Set $\widetilde{A}=\left(\begin{array}{lcr}\pi^+P+G &   K \\
T &  \widetilde{S}    \end{array}\right)$ $\in \mathcal{B}$ , and denote by $p$, $b$ and $s$ the local symbols of $P,G$ and $\widetilde{S}$ respectively.
 Define:
 \begin{align}
{\rm{\widetilde{Wres}}}(\widetilde{A})&=\int_X\int_{\bf \widetilde{ S}}{\rm{tr}}_E\left[p_{-n}(x,\xi)\right]\sigma(\xi)dx \nonumber\\
&+2\pi\int_ {\partial X}\int_{\bf \widetilde{S}'}\left\{{\rm tr}_E\left[({\rm{tr}}b_{-n})(x',\xi')\right]+{\rm{tr}}
_F\left[s_{1-n}(x',\xi')\right]\right\}\sigma(\xi')dx',
\end{align}
where ${\rm{\widetilde{Wres}}}$ denotes the noncommutative residue of an operator in the Boutet de Monvel's algebra.\\
Then~~ a) ${\rm \widetilde{Wres}}([\widetilde{A},B])=0 $, for any
$\widetilde{A},B\in\mathcal{B}$;~~ b) It is the unique continuous trace on
$\mathcal{B}/\mathcal{B}^{-\infty}$.
\end{thm}
\begin{prop}\cite{Wa4}\label{prop1} The following identity holds:
\begin{align}\label{a1}
&1) When~~ p_1+p_2=n,~ Vol_n^{(p_1,p_2)}M=c_0VolM;\\
&2) When~~ p_1+p_2\equiv~ n~ mod ~1,~ Vol_n^{(p_1,p_2)}M=\int_{\partial M}\Psi.
\end{align}
\end{prop}

\section{ The noncommutative residue $\widetilde{{\rm Wres}}[\pi^+(c(X)D^{-1})\circ\pi^+(D^{-(2m-2)})]$ on odd dimensional manifolds with boundary}
\label{section:3}
Firstly we recall the definition of the Dirac operator. Let $M$ be an $n=2m+1$ dimensional oriented
compact spin Riemannian manifold with a Riemannian metric $g^{M}$ and let $\nabla^L$ be the Levi-Civita connection about $g^{M}$. In the fixed orthonormal frame $\{\widetilde{e}_1,\cdots,\widetilde{e}_n\}$, the connection matrix $(\omega_{s,t})$ is defined by
\begin{equation}
\label{a2}
\nabla^L(\widetilde{e}_1,\cdots,\widetilde{e}_n)= (\widetilde{e}_1,\cdots,\widetilde{e}_n)(\omega_{s,t}).
\end{equation}
\indent Let $c(\widetilde{e}_i)$ denotes the Clifford action, which satisfies
\begin{align}
\label{a4}
&c(\widetilde{e}_i)c(\widetilde{e}_j)+c(\widetilde{e}_j)c(\widetilde{e}_i)=-2g^{M}(\widetilde{e}_i,\widetilde{e}_j).
\end{align}
In \cite{Y}, the Dirac operator is given
\begin{align}
\label{a5}
D=\sum^n_{i=1}c(\widetilde{e}_i)\bigg[\widetilde{e}_i-\frac{1}{4}\sum_{s,t}\omega_{s,t}
(\widetilde{e}_i)
c(\widetilde{e}_s)c(\widetilde{e}_t)\bigg].
\end{align}
Set a Clifford action $c(X)$ on $M$ and $X=\sum_{\alpha=1}^na_{\alpha}\widetilde{e}_\alpha=X^T+X_n\partial_{x_n}=\sum_{j=1}^nX_j\partial_j$ is a vector field. We define $\nabla_X^{S(TM)}:=X+\frac{1}{4}\sum_{ij}\langle\nabla_X^L{\widetilde{e}_i},\widetilde{e}_j\rangle c(\widetilde{e}_i)c(\widetilde{e}_j)$, which is a spin connection, where $A(X)=\frac{1}{4}\sum_{ij}\langle\nabla_X^L{\widetilde{e}_i},\widetilde{e}_j\rangle c(\widetilde{e}_i)c(\widetilde{e}_j)$. And let $g^{ij}=g(dx_{i},dx_{j})$, $\xi=\sum_{k}\xi_{j}dx_{j}$ and $\nabla^L_{\partial_{i}}\partial_{j}=\sum_{k}\Gamma_{ij}^{k}\partial_{k}$,  we denote that
\begin{align}
&\sigma_{i}=-\frac{1}{4}\sum_{s,t}\omega_{s,t}
(\widetilde{e}_i)c(\widetilde{e}_i)c(\widetilde{e}_s)c(\widetilde{e}_t)
;~~~\xi^{j}=g^{ij}\xi_{i};~~~~\Gamma^{k}=g^{ij}\Gamma_{ij}^{k};~~~~\sigma^{j}=g^{ij}\sigma_{i}.
\end{align}
Then by \cite{Wa3} and $\sigma(\partial_{x_j})=\sqrt{-1}\xi_j$, we have the following lemmas.
\begin{lem}\label{lem2} The following identities hold:
\begin{align}
&\sigma_1(D)=\sqrt{-1}c(\xi); \nonumber\\
&\sigma_0(D)=-\frac{1}{4}\sum_{i,s,t}\omega_{s,t}(\widetilde{e}_i)c(\widetilde{e}_i)c(\widetilde{e}_s)c(\widetilde{e}_t)\nonumber\\
&\sigma_{0}(\nabla_X^{S(TM)})=A(X);\nonumber\\
&\sigma_{1}(\nabla_X^{S(TM)})=\sqrt{-1}\sum_{j=1}^nX_j\xi_j.\nonumber
\end{align}
\end{lem}
By the composition formula of pseudodifferential operators, we have
\begin{lem}\label{lem3} The following identities hold:
\begin{align}
\label{ab22}
&\sigma_{-1}({D}^{-1})=\frac{\sqrt{-1}c(\xi)}{|\xi|^2};\nonumber\\
&\sigma_{-2}(D^{-2})=|\xi|^{-2};\nonumber\\
&\sigma_{-2}({D}^{-1})=\frac{c(\xi)\sigma_{0}(D)c(\xi)}{|\xi|^4}+\frac{c(\xi)}{|\xi|^6}\sum_jc(dx_j)
\Big[\partial_{x_j}(c(\xi))|\xi|^2-c(\xi)\partial_{x_j}(|\xi|^2)\Big]\nonumber\\
&\sigma_{-3}(D^{-2})=-\sqrt{-1}|\xi|^{-4}\xi_k(\Gamma^k-2\sigma^k)-\sqrt{-1}|\xi|^{-6}2\xi^j\xi_\alpha\xi_\beta\partial_jg^{\alpha\beta}.
\end{align}
\end{lem}
Next, we compute the residue $\widetilde{{\rm Wres}}[\pi^+(c(X){D}^{-1})\circ\pi^+({D}^{-(2m-2)})]$ on odd dimensional oriented
compact spin manifolds with boundary and get a general Kastler-Kalau-Walze
type theorem in this case.

By Proposition \ref{prop1}, we have
\begin{equation}\label{a2}
\widetilde{{\rm Wres}}[\pi^+(c(X){D}^{-1})\circ\pi^+({D}^{-(2m-2)})]=\int_{\partial M}\Psi.
\end{equation}
where
\begin{align}
\label{aqqq15}
\Psi &=\int_{|\xi'|=1}\int^{+\infty}_{-\infty}\sum^{\infty}_{j, k=0}\sum\frac{(-i)^{|\alpha|+j+k+1}}{\alpha!(j+k+1)!}
\times {\rm tr}_{\wedge^*T^*M\bigotimes\mathbb{C}}[\partial^j_{x_n}\partial^\alpha_{\xi'}\partial^k_{\xi_n}\sigma^+_{r}(c(X){D}^{-1})(x',0,\xi',\xi_n)
\nonumber\\
&\times\partial^\alpha_{x'}\partial^{j+1}_{\xi_n}\partial^k_{x_n}\sigma_{l}({D}^{-(2m-2)})(x',0,\xi',\xi_n)]d\xi_n\sigma(\xi')dx',
\end{align}
and the sum is taken over $r+l-k-j-|\alpha|-1=-(2m+1),~~r\leq -1,~~l\leq -(2m-2)$.

 When $n=2m+1$ is odd, then ${\rm tr}_{S(TM)}[{\rm \texttt{id}}]=2^{\frac{n-1}{2}}=2^m$, the sum is taken over $
r+l-k-j-|\alpha|=-2m,~~r\leq -1,~~l\leq -(2m-2),$ then we have the following five cases:
~\\
\noindent  {\bf (case a-I)}~$r=-1,~l=-(2m-2),~k=j=0,~|\alpha|=1$.\\
\noindent By (\ref{aqqq15}), we get
\begin{equation}
\label{ab24}
\Psi_1=-\int_{|\xi'|=1}\int^{+\infty}_{-\infty}\sum_{|\alpha|=1}
 {\rm tr}[\partial^\alpha_{\xi'}\pi^+_{\xi_n}\sigma_{-1}(c(X){D}^{-1})\times
 \partial^\alpha_{x'}\partial_{\xi_n}\sigma_{-(2m-2)}(D^{-(2m-2)})](x_0)d\xi_n\sigma(\xi')dx'.
\end{equation}
By Lemma 2.2 in \cite{Wa3}, for $i<n$, then
\begin{equation}
\label{ab25}
\partial_{x_i}\sigma_{-(2m-2)}(D^{-(2m-2)})(x_0)=\partial_{x_i}{(|\xi|^{(2-2m)})}(x_0)
=\partial_{x_i}(|\xi|^2)^{(1-m)}(x_0)=(1-m) (|\xi|^2)^{-m}
 \partial_{x_i}(|\xi|^2)(x_0)=0,
\end{equation}
\noindent so $\Psi_1=0$.\\
 \noindent  {\bf(case a-II)}~$r=-1,~l=-(2m-2),~k=|\alpha|=0,~j=1$.\\
\noindent By (\ref{aqqq15}), we get
\begin{equation}
\label{ab26}
\Psi_2=-\frac{1}{2}\int_{|\xi'|=1}\int^{+\infty}_{-\infty} {\rm
tr} [\partial_{x_n}\pi^+_{\xi_n}\sigma_{-1}(c(X){D}^{-1})\times
\partial_{\xi_n}^2\sigma_{-(2m-2)}(D^{-(2m-2)})](x_0)d\xi_n\sigma(\xi')dx'.
\end{equation}
\noindent By (3.16) in \cite{Wa6}, we have\\
\begin{align}\label{ab237}
\partial^2_{\xi_n} \sigma_{-(2m-2)}(D^{-(2m-2)})(x_0)
&=\partial^2_{\xi_n}\big((|\xi|^2)^{1-m}\big) (x_0)
=\partial_{\xi_n}\Big( (1-m) (|\xi|^{2})^{-m}  \partial_{\xi_n} (|\xi|^2)\Big) (x_0)\nonumber\\
&=-m(1-m)(|\xi|^{2})^{-m-1} \big(\partial_{\xi_n}|\xi|^2\big)^{2}(x_0)+(1-m)(|\xi|^{2})^{-m} \partial^2_{\xi_n}(|\xi|^2 (x_0))\nonumber\\
&=\Big((4m-2) \xi_n^{2}-2 \Big)(m-1)(1+\xi_n^{2})^{(-m-1)}.
\end{align}
By Lemma \ref{lem3}, we have\\
\begin{align}\label{ab27}
\partial_{x_n}\sigma_{-1}(c(X)D^{-1})(x_0)&=\frac{i\partial_{x_n}(c(X))c(\xi)}{|\xi|^2}+\frac{ic(X)\partial_{x_n}c(\xi')(x_0)}{|\xi|^2}-\frac{ic(X)c(\xi)|\xi'|^2h'(0)}{|\xi|^4}.
\end{align}
By (2.1.1), (2.1.2) in \cite{Wa3} and the Cauchy integral formula, we have
\begin{align}\label{ab29}
\pi^+_{\xi_n}\left[\frac{i\partial_{x_n}(c(X))c(\xi)}{|\xi|^2}\right](x_0)|_{|\xi'|=1}&=i\partial_{x_n}(c(X))\pi^+_{\xi_n}\left[\frac{c(\xi)}{|\xi|^2}\right](x_0)|_{|\xi'|=1}\nonumber\\
&=\partial_{x_n}(c(X))\frac{c(\xi')+ic(dx_n)}{2(\xi_n-i)}.
\end{align}
Similarly, we have
\begin{align}\label{a30}
\pi^+_{\xi_n}\left[\frac{ic(X)\partial_{x_n}c(\xi')(x_0)}{|\xi|^2}\right](x_0)|_{|\xi'|=1}=\frac{c(X)\partial_{x_n}[c(\xi')](x_0)}{2(\xi_n-i)}.
\end{align}
\begin{align}\label{a300}
\pi^+_{\xi_n}\left[\frac{c(X)c(\xi)|\xi'|^2h'(0)}{|\xi|^4}\right](x_0)|_{|\xi'|=1}=-ic(X)\left[\frac{(i\xi_n+2)c(\xi')+ic(dx_n)}{4(\xi_n-i)^2}\right].
\end{align}
Then, we have
\begin{align}\label{ab28}
\partial_{x_n}\pi^+_{\xi_n}\sigma_{-1}(c(X)D^{-1})&=\pi^+_{\xi_n}\partial_{x_n}\sigma_{-1}(c(X)D^{-1})\nonumber\\
&=\partial_{x_n}(c(X))\frac{c(\xi')+ic(dx_n)}{2(\xi_n-i)}+\frac{c(X)\partial_{x_n}[c(\xi')](x_0)}{2(\xi_n-i)}-ic(X)\left[\frac{(i\xi_n+2)c(\xi')+ic(dx_n)}{4(\xi_n-i)^2}\right].
\end{align}
By the relation of the Clifford action and ${\rm tr}{ab}={\rm tr }{ba}$,  we have the equalities:
\begin{align}\label{a49}
&{\rm tr}[c(X)
c(\xi')]=-g(X,\xi'){\rm tr}[\texttt{id}];~~
{\rm tr}[c(X)c(dx_n)]=-X_n{\rm tr}[\texttt{id}];~~{\rm tr}[c(X)
\partial_{x_n}(c(\xi'))]=-\frac{1}{2}h'(0)g(X,\xi'){\rm tr}[\texttt{id}];\nonumber\\
&{\rm tr }[\partial_{x_n}(c(X))c(\xi')]=-\partial_{x_n}(g(X,\xi')){\rm tr}[\texttt{id}]+\frac{1}{2}h'(0)(x_0)g(X,\xi'){\rm tr}[\texttt{id}];
~~{\rm tr}[\partial_{x_n}(c(X))c(dx_n)]=-\partial_{x_n}(X_n){\rm tr}[\texttt{id}].
\end{align}
Then, we have
\begin{align}\label{a3lll3}
&{\rm
tr} [\partial_{x_n}\pi^+_{\xi_n}\sigma_{-1}(c(X)D^{-1})\times
\partial_{\xi_n}^2\sigma_{-(2m-2)}(D^{-(2m-2)})](x_0)\nonumber\\
&=-(m-1)\frac{(4m-2)\xi_n^2-2}{2(\xi_n-i)^{m+2}(\xi_n+i)^{m+1}}\left(\partial_{x_n}(g(X,\xi'))-\frac{1}{2}h'(0)g(X,\xi')(x_0)\right){\rm tr}[\texttt{id}]\nonumber\\
&-i(m-1)\frac{(4m-2)\xi_n^2-2}{2(\xi_n-i)^{m+2}(\xi_n+i)^{m+1}}\partial_{x_n}(X_n){\rm tr}[\texttt{id}]+(m-1)\frac{i\xi_n+2}{4(\xi_n-i)^{m+3}(\xi_n+i)^{m+1}}g(X,\xi'){\rm tr}[\texttt{id}]\nonumber\\
&+(m-1)\frac{i}{4(\xi_n-i)^{m+3}(\xi_n+i)^{m+1}}X_n{\rm tr}[\texttt{id}]-(m-1)\frac{i\xi_n+2}{4(\xi_n-i)^{m+2}(\xi_n+i)^{m+1}}h'(0)g(X,\xi'){\rm tr}[\texttt{id}].
\end{align}
We note that $i<n,~\int_{|\xi'|=1}\xi_{i_{1}}\xi_{i_{2}}\cdots\xi_{i_{2d+1}}\sigma(\xi')=0$,
so $g(X,\xi')$ and $\partial_{x_n}(g(X,\xi'))$ have no contribution for computing ${\bf \Psi_2}$. We have
\begin{align}\label{a35}
\Psi_2&=-\frac{1}{2}\int_{|\xi'|}\int_{-\infty}^{\infty}-i(m-1)\frac{(4m-2)\xi_n^2-2}{2(\xi_n-i)^{m+2}(\xi_n+i)^{m+1}}\partial_{x_n}(X_n){\rm tr}[\texttt{id}]d\xi_n\sigma(\xi')dx'\nonumber\\
&-\frac{1}{2}\int_{|\xi'|}\int_{-\infty}^{\infty}(m-1)\frac{i}{4(\xi_n-i)^{m+3}(\xi_n+i)^{m+1}}X_n{\rm tr}[\texttt{id}]d\xi_n\sigma(\xi')dx'\nonumber\\
&=\frac{i}{2}(m-1) Vol(S^{n-2})\partial_{x_n}(X_n)2^m\int_{\Gamma^+}\frac{(4m-2)\xi_n^2-2}{2(\xi_n-i)^{m+2}(\xi_n+i)^{m+1}}d\xi_ndx'\nonumber\\
&-\frac{i}{8}(m-1)Vol(S^{n-2})X_n2^m \int_{\Gamma^+}\frac{1}{2(\xi_n-i)^{m+2}(\xi_n+i)^{m+1}}d\xi_ndx'\nonumber\\
&=\frac{i}{2}(m-1)Vol(S^{n-2})\partial_{x_n}(X_n)2^m \frac{2\pi i}{(m+1)!}\left[\frac{(4m-2)\xi_n^2-2}{(\xi_n+i)^{m+1}}\right]^{(m+1)}\bigg|_{\xi_n=i}dx'\nonumber\\
&-\frac{i}{8}(m-1)Vol(S^{n-2})X_n2^m \frac{2\pi i}{(m+2)!}\left[\frac{1}{(\xi_n+i)^{m+1}}\right]^{(m+2)}\bigg|_{\xi_n=i}dx'\nonumber\\
  &:=(m-1)\bigg(\partial_{x_n}(X_n)\frac{\pi}{(m+1)!}A_0+X_n\frac{\pi i}{2(m+2)!}A_1\bigg) Vol(S^{n-2})2^mdx',
\end{align}
where let $C_N^K=\frac{N!}{K!(N-K)!}$ and $A_N^K=\frac{N!}{(N-K)!}$, we have
\begin{align}
A_0&=\left[\frac{(4m-2)\xi_n^2-2}{(\xi_n+i)^{m+1}}\right]^{(m+1)}\bigg|_{\xi_n=i}\nonumber\\
 &=-i^{-2(m+1)}2^{-(2m+1)}\bigg((4m-2)C_{-(m+1)}^{m-1}+(4m-2)C_{-(m+1)}^{m}+mC_{-(m+1)}^{m+1}\bigg)(m+1)!;\nonumber\\
A_1&=\left[\frac{1}{(\xi_n+i)^{m+1}}\right]^{(m+2)}\bigg|_{\xi_n=i}=(2i)^{-2m-3}A_{-m-1}^{m+2}.
\end{align}
\noindent  {\bf (case a-III)}~$r=-1,~l=-(2m-2),~j=|\alpha|=0,~k=1$.\\
\noindent By (\ref{aqqq15}), we get
\begin{align}\label{a36}
\Psi_3&=-\frac{1}{2}\int_{|\xi'|=1}\int^{+\infty}_{-\infty}
{\rm tr} [\partial_{\xi_n}\pi^+_{\xi_n}\sigma_{-1}(c(X)D^{-1})\times
\partial_{\xi_n}\partial_{x_n}\sigma_{-(2m-2)}(D^{-(2m-2)})](x_0)d\xi_n\sigma(\xi')dx'\nonumber\\
&=\frac{1}{2}\int_{|\xi'|=1}\int^{+\infty}_{-\infty}
{\rm tr} [\partial_{\xi_n}^2\pi^+_{\xi_n}\sigma_{-1}(c(X)D^{-1})\times
\partial_{x_n}\sigma_{-(2m-2)}(D^{-(2m-2)})](x_0)d\xi_n\sigma(\xi')dx'.
\end{align}
By (3.21) in \cite{Wa6}, we have
\begin{align}\label{a37}
\partial_{x_n} \big(\sigma_{-(2m-2)}(D^{-(2m-2)})\big)(x_0)
=\partial_{x_n}\big((|\xi|^2)^{1-m}\big) (x_0)
=h'(0)(1-m)(1+\xi_{n}^{2})^{-m}.
\end{align}
Then by Lemma \ref{lem3}, we get
\begin{align}\label{ammmmm}
&\partial_{\xi_n}^2\pi^+_{\xi_n}\sigma_{-1}(c(X)D^{-1})(x_0)|_{|\xi'|=1}=c(X)\frac{c(\xi')+ic(dx_n)}{(\xi_n-i)^3}.
\end{align}
Moreover
\begin{align}\label{wwwl3}
&{\rm tr} [\partial_{\xi_n}^2\pi^+_{\xi_n}\sigma_{-1}(c(X)D^{-1})\times
\partial_{x_n}\sigma_{-(2m-2)}(D^{-(2m-2)})](x_0)\nonumber\\
&=-(1-m)h'(0)\frac{1}{(\xi_n-i)^{m+3}(\xi_n+i)^{m}}g(X,\xi'){\rm tr}[\texttt{id}]-(1-m)h'(0)\frac{i}{(\xi_n-i)^{m+3}(\xi_n+i)^{m}}X_n{\rm tr}[\texttt{id}].
\end{align}
Next, we perform the corresponding integral calculation on the above results. When we omit $g(X,\xi')$ and $\partial_{x_n}(g(X,\xi'))$ that have no contribution for computing ${\bf \Psi_3}$, we obtain
\begin{align}\label{a3oooo5}
\Psi_3&= \frac{1}{2}\int_{|\xi'|}\int_{-\infty}^{\infty}-(1-m)h'(0)\frac{i}{(\xi_n-i)^{m+3}(\xi_n+i)^{m}}X_n{\rm tr}[\texttt{id}]d\xi_n\sigma(\xi')dx'\nonumber\\
&=- \frac{1}{2}Vol(S^{n-2})(1-m)X_nh'(0)2^m\int_{\Gamma^+}\frac{i}
   {(\xi_n-i)^{m+3}(\xi_n+i)^{m}}d\xi_ndx'\nonumber\\
&= -\frac{1}{2}Vol(S^{n-2})(1-m)X_nh'(0) 2^m\frac{2\pi i}{(m+2)!}
\left[\frac{i}{(\xi_n+i)^m}\right]^{(m+2)}\bigg|_{\xi_n=i}dx'\nonumber\\
  &:=(m-1)Vol(S^{n-2})X_nh'(0)2^m\frac{\pi i}{(m+2)!}B_0dx',
\end{align}
where
\begin{align}
B_0=\left[\frac{i}{(\xi_n+i)^{m}}\right]
 ^{(m+2)}\bigg|_{\xi_n=i}=-i^{-2m-3}2^{-2m-2}A_{-m}^{m+2}.\nonumber\\
\end{align}
\noindent  {\bf (case a-IV)}~$r=-1,~l=-(2m-1),~k=j=|\alpha|=0$.\\
\noindent By (\ref{aqqq15}), we get
\begin{align}\label{a42}
\Psi_4&=-i\int_{|\xi'|=1}\int^{+\infty}_{-\infty}{\rm tr} [\pi^+_{\xi_n}\sigma_{-1}(c(X)D^{-1})\times
\partial_{\xi_n}\sigma_{-(2m-1)}(D^{-(2m-2)})](x_0)d\xi_n\sigma(\xi')dx'\nonumber\\
&=i\int_{|\xi'|=1}\int^{+\infty}_{-\infty}{\rm tr} [\partial_{\xi_n}\pi^+_{\xi_n}\sigma_{-1}(c(X)D^{-1})\times
\sigma_{-(2m-1)}(D^{-(2m-2)})](x_0)d\xi_n\sigma(\xi')dx'.
\end{align}
By (3.30) in \cite{Wa6}, we have
\begin{align}\label{a43}
\sigma_{-(2m-1)}(D^{-2m+2})
&=(m-1)(1+\xi_{n}^{2})^{(-m+2)}
   \Big[\frac{-i}{(1+\xi_n^2)^2}\times \frac{2m+1}{2}h'(0)\xi_n-\frac{2ih'(0)\xi_n}{(1+\xi_n^2)^3} \Big]\nonumber\\
   &+\sqrt{-1} h'(0)(-m^2+3m-2)  \xi_{n} (1+\xi_{n}^{2})^{(-m-1)}.
\end{align}
 By Lemma \ref{lem3}, we have
\begin{align}\label{a45}
\partial_{\xi_n}\pi^+_{\xi_n}\sigma_{-1}(c(X)D^{-1})&=-c(X)\frac{c(\xi')+ic(dx_n)}{2(\xi_n-i)^2}.
\end{align}
Then, we have
\begin{align}\label{a3lll3}
&{\rm
tr} [\partial_{\xi_n}\pi^+_{\xi_n}\sigma_{-1}(c(X)D^{-1})\times
\sigma_{-(2m-2)}(D^{-(2m-2)})](x_0)\nonumber\\
&=-(2m^2-m-1)h'(0)\frac{i\xi_n}{4(\xi_n-i)^{m+2}(\xi_n+i)^{m}}g(X,\xi'){\rm tr}[\texttt{id}]\nonumber\\
&+(2m^2-m-1)h'(0)\frac{\xi_n}{4(\xi_n-i)^{m+2}(\xi_n+i)^{m}}X_n{\rm tr}[\texttt{id}]\nonumber\\
&-(m^2-2m+1)h'(0)\frac{i\xi_n}{(\xi_n-i)^{m+3}(\xi_n+i)^{m+1}}g(X,\xi'){\rm tr}[\texttt{id}]\nonumber\\
&+(m^2-2m+1)h'(0)\frac{\xi_n}{(\xi_n-i)^{m+3}(\xi_n+i)^{m+1}}X_n{\rm tr}[\texttt{id}].
\end{align}
Similarly, we omit $g(X,\xi')$ that has no contribution for computing ${\bf \Psi_4}$. Then, we get
\begin{align}\label{a39}
\Psi_4&=i\int_{|\xi'|}\int^{+\infty}_{-\infty}(2m^2-m-1)h'(0)\frac{\xi_n}{4(\xi_n-i)^{m+2}(\xi_n+i)^{m}}X_n{\rm tr}[\texttt{id}]d\xi_n\sigma(\xi')dx' \nonumber\\
&+i\int_{|\xi'|}\int^{+\infty}_{-\infty}(m^2-2m+1)h'(0)\frac{\xi_n}{(\xi_n-i)^{m+3}(\xi_n+i)^{m+1}}X_n{\rm tr}[\texttt{id}]d\xi_n\sigma(\xi')dx' \nonumber\\
&=i(2m^2-m-1)Vol(S^{n-2})X_nh'(0)2^m\frac{1}{4}\int_{\Gamma^+}\frac{\xi_n}{4(\xi_n-i)^{m+2}(\xi_n+i)^{m}}d\xi_ndx'\nonumber\\
&+i(m^2-2m+1)Vol(S^{n-2})X_nh'(0)2^m\int_{\Gamma^+}\frac{\xi_n}{(\xi_n-i)^{m+3}(\xi_n+i)^{m+1}}d\xi_ndx'\nonumber\\
&=i(2m^2-m-1)Vol(S^{n-2})X_nh'(0)2^m\frac{1}{8}\frac{2\pi i}{(m+1)!}\left[\frac{\xi_n}{(\xi_n+i)^{m}}\right]^{(m+1)}\bigg|_{\xi_n=i}dx'\nonumber\\
&+i(m^2-2m+1)Vol(S^{n-2})X_nh'(0)2^m\frac{2\pi i}{(m+2)!}\left[\frac{\xi_n}{(\xi_n+i)^{m+1}}\right]^{(m+2)}\bigg|_{\xi_n=i}dx'\nonumber\\
   &:=-\bigg(\frac{(2m^2-m-1)\pi }{4(m+1)!}C_{0}+\frac{2(m^2-2m+1)\pi }{(m+2)!}C_1\bigg)Vol(S^{n-2})X_nh'(0)2^m dx',
\end{align}
where
\begin{align}
C_0&=\left[\frac{\xi_n}{(\xi_n+i)^{m}}\right]^{(m+1)}\bigg|_{\xi_n=i}=-i^{-2m-2}2^{-2m-1}\bigg(2C_{-m}^{m}+C_{-m}^{m+1}\bigg)(m+1)!;\nonumber\\
 C_1&=\left[\frac{\xi_n}{(\xi_n+i)^{m+1}}\right]^{(m+2)}\bigg|_{\xi_n=i}dx'=-i^{-2m-4}2^{-2m-3}\bigg(2C_{-m-1}^{m+1}+C_{-m-1}^{m+2}\bigg)(m+2)!.
\end{align}
\noindent {\bf (case a-V)}~$r=-2,~\ell=-(2m-2),~k=j=|\alpha|=0$.\\
By (\ref{aqqq15}), we get
\begin{align}\label{a61}
\Psi_5=-i\int_{|\xi'|=1}\int^{+\infty}_{-\infty}{\rm tr} [\pi^+_{\xi_n}\sigma_{-2}(c(X)D^{-1})\times
\partial_{\xi_n}\sigma_{-(2m-2)}(D^{-(2m-2)})](x_0)d\xi_n\sigma(\xi')dx'.
\end{align}
By (3.33) in \cite{Wa6}, we have
\begin{equation}\label{a62}
\partial_{\xi_n}\sigma_{-(2m-2)}(D^{-(2m-2)})(x_0)
=\partial_{\xi_n}((|\xi|^2)^{1-m})(x_0)=2(1-m)\xi_n(1+\xi_n^2)^{-m}.
\end{equation}
By Lemma \ref{lem3}, we have
\begin{align}\label{a43}
\sigma_{-2}(c(X)D^{-1})(x_0)=c(X)\left\{\frac{c(\xi)\sigma_{0}(D)c(\xi)}{|\xi|^4}+\frac{c(\xi)}{|\xi|^6}\sum_jc(dx_j)
\Big[\partial_{x_j}(c(\xi))|\xi|^2-c(\xi)\partial_{x_j}(|\xi|^2)\Big]\right\},
\end{align}
where
\begin{align}\label{a44}
\sigma_{0}(D)(x_0)&=-\frac{1}{4}\sum_{s,t,i}\omega_{s,t}(\widetilde{e}_i)
(x_{0})c(\widetilde{e}_i)c(\widetilde{e}_s)c(\widetilde{e}_t).
\end{align}
We denote
\begin{align}\label{a45}
H(x_0)&=-\frac{1}{4}\sum_{s,t,i}\omega_{s,t}(\widetilde{e}_i)
(x_{0})c(\widetilde{e}_i)c(\widetilde{e}_s)c(\widetilde{e}_t),
\end{align}
where $H(x_0)=c_0c(dx_n)$ and $c_0=-\frac{3}{4}h'(0)$.\\
Moreover
\begin{align}\label{a46}
&\pi^+_{\xi_n}\sigma_{-2}(c(X)D^{-1}(x_0))|_{|\xi'|=1}\nonumber\\
&=
\pi^+_{\xi_n}\Big[\frac{c(X)c(\xi)H(x_0)c(\xi)+c(X)c(\xi)c(dx_n)\partial_{x_n}[c(\xi')](x_0)}{(1+\xi_n^2)^2}-h'(0)\frac{c(X)c(\xi)c(dx_n)c(\xi)}{(1+\xi_n^{2})^3}\Big].
\end{align}
By computations, we have
\begin{align}\label{a52}
\pi^+_{\xi_n}\Big[\frac{c(X)c(\xi)H(x_0)c(\xi)+c(X)c(\xi)c(dx_n)\partial_{x_n}[c(\xi')](x_0)}{(1+\xi_n^2)^2}\Big]-h'(0)\pi^+_{\xi_n}\Big[c(X)\frac{c(\xi)c(dx_n)c(\xi)}{(1+\xi_n)^3}\Big]:= E^1-E^2,
\end{align}
where
\begin{align}\label{a53}
E^1&=\frac{-1}{4(\xi_n-i)^2}[(2+i\xi_n)c(X)c(\xi')H(x_0)c(\xi')+i\xi_nc(X)c(dx_n)H(x_0)c(dx_n)\nonumber\\
&+(2+i\xi_n)c(X)c(\xi')c(dx_n)\partial_{x_n}c(\xi')+ic(X)c(dx_n)H(x_0)c(\xi')
+ic(X)c(\xi')H(x_0)c(dx_n)-ic(X)\partial_{x_n}c(\xi')]
\end{align}
and
\begin{align}\label{a54}
E^2&=\frac{h'(0)}{2}c(X)\left[\frac{c(dx_n)}{4i(\xi_n-i)}+\frac{c(dx_n)-ic(\xi')}{8(\xi_n-i)^2}
+\frac{3\xi_n-7i}{8(\xi_n-i)^3}[ic(\xi')-c(dx_n)]\right].
\end{align}
Then, we have
\begin{align}\label{ajjl3}
&{\rm
tr} [\pi^+_{\xi_n}\sigma_{-2}(c(X)D^{-1})\times
\partial_{\xi_n}\sigma_{-(2m-2)}(D^{-(2m-2)})](x_0)\nonumber\\
&=(1-m)h'(0)\frac{\xi_n^2+5i\xi_n}{4(\xi_n-i)^{m+3}(\xi_n+i)^{m}}X_n{\rm tr}[\texttt{id}]+(1-m)h'(0)\frac{2i\xi_n^2+3\xi_n}{(\xi_n-i)^{m+3}(\xi_n+i)^{m}}g(X,\xi'){\rm tr}[\texttt{id}].
\end{align}
We omit $g(X,\xi')$ that has no contribution for computing ${\bf \Psi_5}$. Then, we obtain
\begin{align}\label{a66}
\Psi_5&=-i\int_{|\xi'|=1}\int^{+\infty}_{-\infty}(1-m)h'(0)\frac{\xi_n^2+5i\xi_n}{4(\xi_n-i)^{m+3}(\xi_n+i)^{m}}X_n{\rm tr}[\texttt{id}]d\xi_n\sigma(\xi')dx'\nonumber\\
&=-i(1-m)Vol(S^{n-2})X_nh'(0)2^m\frac{1}{4}\int_{\Gamma^+}\frac{\xi_n^2+5i\xi_n}{(\xi_n-i)^{m+3}(\xi_n+i)^{m}}dx'\nonumber\\
&=-i(1-m)Vol(S^{n-2})X_nh'(0)2^m\frac{1}{4}\frac{2\pi i}{(m+2)!}\left[\frac{\xi_n^2+5i\xi_n}{(\xi_n+i)^{m}}\right]^{(m+2)}\bigg|_{\xi_n=i}dx'\nonumber\\
&:=-(m-1)Vol(S^{n-2})X_nh'(0)2^m\frac{\pi}{2(m+2)!}D_0dx',
\end{align}
where
\begin{align}
D_0&=\left[\frac{\xi_n^2+5i\xi_n}{(\xi_n+i)^{m}}\right]^{(m+2)}\bigg|_{\xi_n=i}=-i^{-2m-3}2^{-2m-1}\bigg(2C_{-m}^{m}+7C_{-m}^{m+1}+3C_{-m}^{m+2}\bigg)(m+2)!.
\end{align}
Now $\Psi$ is the sum of the {\bf  (case a-I)}-{\bf (case a-V)}. Therefore, we get
\begin{align}\label{a795}
\Psi&=\sum_{i=1}^5\Psi_i\nonumber\\
&=\bigg\{(m-1)\bigg(\partial_{x_n}(X_n)\frac{\pi}{(m+1)!}A_0+\frac{\pi i}{2(m+2)!}A_1+X_nh'(0)\frac{\pi i}{(m+2)!}B_0\bigg) -X_nh'(0)\frac{(2m^2-m-1)\pi}{4(m+1)!}C_{0}\nonumber\\
&-X_nh'(0)\frac{2(m^2-2m+1)\pi }{(m+2)!}C_1
-(m-1)X_nh'(0)\frac{\pi }{2(m+2)!}D_0\bigg\}2^m Vol(S^{n-2})dx'.\nonumber\\
\end{align}
Then, by (\ref{a2}) and (\ref{a795}), we obtain following theorem
\begin{thm}\label{thmb1}
Let $M$ be an $n=2m+1$ dimensional oriented
compact spin manifold with boundary $\partial M$, then we get the following equality:
\begin{align}
\label{ab2773}
&\widetilde{{\rm Wres}}[\pi^+(c(X)D^{-1})\circ\pi^+(D^{-(2m-2)})]\nonumber\\
&=\int_{\partial M} \bigg(\partial_{x_n}(X_n)\frac{\pi}{(m+1)!}A_0+X_n\frac{\pi i}{2(m+2)!}A_1+X_nh'(0)\frac{\pi i}{(m+2)!}B_0\bigg) -X_nh'(0)\frac{(2m^2-m-1)\pi}{4(m+1)!}C_{0}\nonumber\\
&-X_nh'(0)\frac{2(m^2-2m+1)\pi }{(m+2)!}C_1
-(m-1)X_nh'(0)\frac{\pi }{2(m+2)!}D_0\bigg\}Vol(S^{n-2})2^md{\rm Vol_{M}}.
\end{align}
\end{thm}
\section{The noncommutative residue $\widetilde{{\rm Wres}}[\pi^+(\nabla_X^{S(TM)}D^{-1})\circ\pi^+(D^{-(2m-1)})]$ on odd dimensional manifolds with boundary}
\label{section:4}
In this section, we compute the noncommutative residue $\widetilde{{\rm Wres}}[\pi^+(\nabla_X^{S(TM)}D^{-1})\circ\pi^+(D^{-(2m-1)})]$ on $2m+1$ dimensional oriented
compact spin manifolds with boundary.

Similar to \cite{Wa3}, by Proposition \ref{prop1}, we can compute the noncommutative residue
\begin{align}
\label{c1}
&\widetilde{{\rm Wres}}[\pi^+(\nabla_X^{S(TM)}D^{-1})\circ\pi^+(D^{-(2m-1)})]=\int_{\partial M}\widetilde{\Psi},
\end{align}
where
\begin{align}
\label{c2}
 \widetilde{\Psi}&=\int_{|\xi'|=1}\int^{+\infty}_{-\infty}\sum^{\infty}_{j, k=0}\sum\frac{(-i)^{|\alpha|+j+k+1}}{\alpha!(j+k+1)!}
\times {\rm tr}_{\wedge^*T^*M\bigotimes\mathbb{C}}[\partial^j_{x_n}\partial^\alpha_{\xi'}\partial^k_{\xi_n}\sigma^+_{r}(\nabla_X^{S(TM)}D^{-1})(x',0,\xi',\xi_n)
\nonumber\\
&\times\partial^\alpha_{x'}\partial^{j+1}_{\xi_n}\partial^k_{x_n}\sigma_{l}(D^{-(2m-1)})(x',0,\xi',\xi_n)]d\xi_n\sigma(\xi')dx',
\end{align}
and the sum is taken over $r+l-k-j-|\alpha|-1=-(2m+1),~~r\leq 0,~~l\leq -(2m-1)$.

 Therefore, we  need to compute $\int_{\partial M} \widetilde{\Psi}$. When sum is taken over $r+l-k-j-|\alpha|=-2m,~~r\leq 0,~~l\leq -(2m-1),$  we have the following five cases:

\noindent  {\bf(case b-I)}~$r=0,~l=-(2m-1),~k=j=0,~|\alpha|=1$.\\
\noindent By (\ref{c2}), we get
\begin{equation}
\label{c4}
\widetilde{\Psi}_1=-\int_{|\xi'|=1}\int^{+\infty}_{-\infty}\sum_{|\alpha|=1}
{\rm tr}[\partial^\alpha_{\xi'}\pi^+_{\xi_n}\sigma_{0}(\nabla_X^{S(TM)}D^{-1})\times
 \partial^\alpha_{x'}\partial_{\xi_n}\sigma_{-(2m-1)}(D^{-(2m-1)})](x_0)d\xi_n\sigma(\xi')dx'.
\end{equation}
By Lemma 2.2 in \cite{Wa3}, for $i<n$, then
\begin{align}
\label{c5}
\partial_{x_i}\sigma_{-(2m-1)}(D^{-(2m-1)})(x_0)&=\partial_{x_i}{(ic(\xi)|\xi|^{-2m})}(x_0)\nonumber\\
&=i\partial_{x_i}c(\xi)(x_0)|\xi|^{-2m}+ic(\xi)\partial_{x_i}(|\xi|^{-2m})(x_0)\nonumber\\
 &=0,
\end{align}
\noindent so $\widetilde{\Psi}_1=0$.\\
 \noindent  {\bf (case b-II)}~$r=0,~l=-(2m-1),~k=|\alpha|=0,~j=1$.\\
\noindent By (\ref{c2}), we get
\begin{align}
\label{c6}
\widetilde{\Psi}_2&=-\frac{1}{2}\int_{|\xi'|=1}\int^{+\infty}_{-\infty} {\rm
tr} [\partial_{x_n}\pi^+_{\xi_n}\sigma_{0}(\nabla_X^{S(TM)}D^{-1})\times
\partial_{\xi_n}^2\sigma_{-(2m-1)}(D^{-(2m-1)})](x_0)d\xi_n\sigma(\xi')dx'\nonumber\\
&=-\frac{1}{2}\int_{|\xi'|=1}\int^{+\infty}_{-\infty} {\rm
tr} [\partial_{\xi_n}^2\partial_{x_n}\pi^+_{\xi_n}\sigma_{0}(\nabla_X^{S(TM)}D^{-1})\times
\sigma_{-(2m-1)}(D^{-(2m-1)})](x_0)d\xi_n\sigma(\xi')dx'.
\end{align}
\noindent By (3.11) in \cite{Wa7}, we have\\
\begin{eqnarray}\label{c7}
\sigma_{-(2m-1)}(D^{-(2m-1)})
&=&\frac{\sqrt{-1}[c(\xi')+\xi_nc(\mathrm{d}x_n)]}{(1+\xi_n^2)^{m}}.
\end{eqnarray}
By Lemma \ref{lem2} and Lemma \ref{lem3}, we have
\begin{align}\label{c8}
&\partial_{x_n}\sigma_{0}(\nabla_X^{S(TM)}D^{-1})\nonumber\\
&=\partial_{x_n}\bigg(-\sum_{j,l=1}^nX_j\xi_j\frac{c(\xi)}{|\xi|^{2}}\bigg)\nonumber\\
&=-\sum_{j=1}^{n-1}\xi_j\left[\frac{\partial_{x_n}(X_jc(\xi))}{|\xi|^{2}}-\frac{h'(0)|\xi'|^2X_jc(\xi)}{|\xi|^{4}}+\frac{\xi_n\partial_{x_n}(X_jc(dx_n))}{|\xi|^{2}}-\frac{\xi_nh'(0)|\xi'|^2X_jc(dx_n)}{|\xi|^{4}}\right]\nonumber\\
&-\frac{\xi_n\partial_{x_n}(X_nc(\xi))}{|\xi|^{2}}+\frac{\xi_nh'(0)|\xi'|^2X_nc(\xi)}{|\xi|^{4}}-\frac{\xi_n^2\partial_{x_n}(X_nc(dx_n))}{|\xi|^{2}}-\frac{\xi_n^2h'(0)|\xi'|^2X_nc(dx_n)}{|\xi|^{4}}.
\end{align}
 Then, we have
\begin{align}\label{c9}
&\partial_{x_n}\pi^+_{\xi_n}\sigma_{0}(\nabla_X^{S(TM)}D^{-1})\nonumber\\
&=\pi^+_{\xi_n}\partial_{x_n}\sigma_{0}(\nabla_X^{S(TM)}D^{-1})\nonumber\\
&=\frac{i}{2(\xi_n-i)}\sum_{j=1}^{n-1}\xi_j\partial_{x_n}(X_jc(\xi'))+\frac{2+i\xi_n}{4(\xi_n-i)^2}\sum_{j=1}^{n-1}\xi_jh'(0)|\xi'|^2X_jc(\xi')+\frac{1}{2(\xi_n-i)}\sum_{j=1}^{n-1}\xi_j\partial_{x_n}(X_jc(dx_n))\nonumber\\
&+\frac{i}{4(\xi_n-i)^2}\sum_{j=1}^{n-1}\xi_jh'(0)|\xi'|^2X_jc(dx_n)-\frac{i}{2(\xi_n-i)}\partial_{x_n}(X_jc(dx_n))+\frac{i}{4(\xi_n-i)^2}h'(0)|\xi'|^2X_nc(dx_n)\nonumber\\
&-\frac{1}{2(\xi_n-i)}\partial_{x_n}(X_nc(\xi'))-\frac{i}{4(\xi_n-i)^2}h'(0)|\xi'|^2X_nc(\xi').
\end{align}
By further calculation, we have
\begin{align}\label{555}
&\partial_{\xi_n}^2\partial_{x_n}\pi^+_{\xi_n}\sigma_{0}(\nabla_X^{S(TM)}D^{-1})\nonumber\\
&=\frac{i}{(\xi_n-i)^3}\sum_{j=1}^{n-1}\xi_j\partial_{x_n}(X_jc(\xi'))+\frac{4+i\xi_n}{2(\xi_n-i)^4}\sum_{j=1}^{n-1}\xi_jh'(0)|\xi'|^2X_jc(\xi')+\frac{1}{(\xi_n-i)^3}\sum_{j=1}^{n-1}\xi_j\partial_{x_n}(X_jc(dx_n))\nonumber\\
&+\frac{3i}{2(\xi_n-i)^4}\sum_{j=1}^{n-1}\xi_jh'(0)|\xi'|^2X_jc(dx_n)-\frac{i}{(\xi_n-i)^3}\partial_{x_n}(X_jc(dx_n))+\frac{2i\xi_n-4}{4(\xi_n-i)^4}h'(0)|\xi'|^2X_nc(dx_n)\nonumber\\
&-\frac{1}{(\xi_n-i)^3}\partial_{x_n}(X_nc(\xi'))-\frac{3i}{2(\xi_n-i)^2}h'(0)|\xi'|^2X_nc(\xi').
\end{align}
By the relation of the Clifford action and ${\rm tr}{ab}={\rm tr }{ba}$,  we have the equalities:
\begin{align}\label{a49}
&{\rm tr}[\partial_{x_n}(X_jc(\xi'))c(dx_n)]={\rm tr}[\partial_{x_n}(X_jc(dx_n))c(\xi')]={\rm tr}[c(\xi')c(dx_n)]=0;\nonumber\\
&{\rm tr}[c(dx_n)c(dx_n)]={\rm tr}[c(\xi')c(\xi')]=-{\rm tr}[\texttt{id}];~~~{\rm tr}[\partial_{x_n}(X_nc(dx_n))c(dx_n)]=-\partial_{x_n}(X_n){\rm tr}[\texttt{id}];\nonumber\\
&{\rm tr}[\partial_{x_n}(X_nc(\xi'))c(\xi')]=-[\partial_{x_n}(X_n)+\frac{1}{2}h'(0)X_n]{\rm tr}[\texttt{id}].
\end{align}
We omit some items that have no contribution for computing ${\bf \widetilde{ \Psi}_2}$.
Then, we have
\begin{align}\label{c3lll3}
&{\rm
tr} [\partial_{\xi_n}^2\partial_{x_n}\pi^+_{\xi_n}\sigma_{0}(\nabla_X^{S(TM)}D^{-1})\times
\sigma_{-(2m-1)}(D^{-(2m-1)})](x_0)\nonumber\\
&=\frac{i-\xi_n}{2(\xi_n-i)^{m+3}(\xi_n+i)^{m}}\partial_{x_n}(X_n){\rm tr}[\texttt{id}]+\frac{-2\xi_n^2+9i\xi_n+4}{2(\xi_n-i)^{m+4}(\xi_n+i)^{m}}h'(0)X_n{\rm tr}[\texttt{id}].
\end{align}
Therefore, we get
\begin{align}\label{c11}
\widetilde{\Psi}_2&=-\frac{1}{2}\int_{|\xi'|}\int_{-\infty}^{\infty}\frac{i-\xi_n}{2(\xi_n-i)^{m+3}(\xi_n+i)^{m}}\partial_{x_n}(X_n){\rm tr}[\texttt{id}]d\xi_n\sigma(\xi')dx'\nonumber\\
&-\frac{1}{2}\int_{|\xi'|}\int_{-\infty}^{\infty}\frac{-2\xi_n^2+9i\xi_n+4}{2(\xi_n-i)^{m+4}(\xi_n+i)^{m}}X_nh'(0){\rm tr}[\texttt{id}]d\xi_n\sigma(\xi')dx'\nonumber\\
 &=-\frac{1}{2}Vol(S^{n-2})\partial_{x_n}(X_n)2^m\int_{\Gamma^+}\frac{i-\xi_n}{2(\xi_n-i)^{m+3}(\xi_n+i)^{m}}d\xi_ndx'\nonumber\\
 &-\frac{1}{2}Vol(S^{n-2})X_nh'(0)2^m\int_{\Gamma^+}\frac{-2\xi_n^2+9i\xi_n+4}{2(\xi_n-i)^{m+4}(\xi_n+i)^{m}}d\xi_ndx'\nonumber\\
 &=-\frac{1}{4}Vol(S^{n-2})\partial_{x_n}(X_n)2^m\frac{2\pi i}{(m+2)!}\left[\frac{i-\xi_n}{(\xi_n+i)^{m}}\right]^{(m+2)}\bigg|_{\xi_n=i}dx'\nonumber\\
 &-\frac{1}{4}Vol(S^{n-2})X_nh'(0)2^m\frac{2\pi i}{(m+3)!}\left[\frac{-2\xi_n^2+9i\xi_n+4}{(\xi_n+i)^{m}}\right]^{(m+3)}\bigg|_{\xi_n=i}dx'\nonumber\\
 &:=-\bigg(\partial_{x_n}(X_n)\frac{\pi i}{2(m+2)!}E_0+h'(0)X_n\frac{\pi i}{2(m+3)!}E_1\bigg)Vol(S^{n-2})2^mdx',
\end{align}
where
\begin{align}
E_0&=\left[\frac{i-\xi_n}{(\xi_n+i)^{m}}\right]^{(m+2)}\bigg|_{\xi_n=i}=i^{-2m+1}2^{-2m-1}C_{-m}^{m+1}(m+2)!;\nonumber\\
E_1&=\left[\frac{-2\xi_n^2+9i\xi_n+4}{(\xi_n+i)^{m}}\right]^{(m+3)}\bigg|_{\xi_n=i}=(2i)^{-2m-3}\bigg(8iC_{-m}^{m+1}-(18-8i)C_{-m}^{m+2}-(5-2i)C_{-m}^{m+3}\bigg)(m+3)!.
\end{align}
\noindent  {\bf (case b-III)}~$r=0,~l=-(2m-1),~j=|\alpha|=0,~k=1$.\\
\noindent By (\ref{c2}), we get
\begin{align}\label{c12}
\widetilde{\Psi}_3&=-\frac{1}{2}\int_{|\xi'|=1}\int^{+\infty}_{-\infty}
{\rm tr} [\partial_{\xi_n}\pi^+_{\xi_n}\sigma_{0}(\nabla_X^{S(TM)}D^{-1})\times
\partial_{\xi_n}\partial_{x_n}\sigma_{-(2m-1)}(D^{-(2m-1)})](x_0)d\xi_n\sigma(\xi')dx'\nonumber\\
&=\frac{1}{2}\int_{|\xi'|=1}\int^{+\infty}_{-\infty}
{\rm tr} [\partial_{\xi_n}^2\pi^+_{\xi_n}\sigma_{0}(\nabla_X^{S(TM)}D^{-1})\times
\partial_{x_n}\sigma_{-(2m-1)}(D^{-(2m-1)})](x_0)d\xi_n\sigma(\xi')dx'.
\end{align}
\noindent By (3.17) in \cite{Wa7}, we have
\begin{eqnarray}\label{c13}
\partial_{x_n} \big(\sigma_{-(2m-1)}(D^{-(2m-1)})\big)(x_0)
=\frac{i\partial_{x_n}c(\xi')(x_0)}{(1+\xi_n^2)^m}-\frac{2mih'(0)c(\xi)}{2(1+\xi_n^2)^{m+1}}.
\end{eqnarray}
By Lemma \ref{lem2} and Lemma \ref{lem3}, we have
\begin{align}\label{c14}
\partial_{\xi_n}^2\pi^+_{\xi_n}\sigma_{0}(\nabla_X^{S(TM)}D^{-1})&=\frac{i}{(\xi_n-i)^3}\sum_{j=1}^{n-1}X_j\xi_jc(\xi')-\frac{1}{(\xi_n-i)^3}\sum_{j=1}^{n-1}X_j\xi_jc(dx_n)\nonumber\\
&-\frac{1}{(\xi_n-i)^3}X_nc(\xi')-\frac{i}{(\xi_n-i)^3}X_nc(dx_n).
\end{align}
Then by (\ref{a49}), we have
\begin{align}\label{c16}
{\rm tr} [\partial_{\xi_n}^2\pi^+_{\xi_n}\sigma_{0}(\nabla_X^{S(TM)}D^{-1})\times
\partial_{x_n}\sigma_{-(2m-1)}(D^{-(2m-1)})](x_0)= \frac{i\xi_n^2+2m(1-i)\xi_n+2m+i}{2(\xi_n-i)^{m+4}(\xi_n+i)^{m+1}}X_nh'(0){\rm tr}[\texttt{id}].
\end{align}
Therefore, we get
\begin{align}\label{c18}
\widetilde{\Psi}_3&=\frac{1}{2}\int_{|\xi'|=1}\int^{+\infty}_{-\infty}
{\rm tr} [\partial_{\xi_n}^2\pi^+_{\xi_n}\sigma_{0}(\nabla_X^{S(TM)}D^{-1})\times
\partial_{x_n}\sigma_{-(2m-1)}(D^{-(2m-1)})](x_0)d\xi_n\sigma(\xi')dx'\nonumber\\
&=\frac{1}{2}\int_{|\xi'|=1}\int^{+\infty}_{-\infty}
\frac{i\xi_n^2+2m(1-i)\xi_n+2m+i}{2(\xi_n-i)^4(\xi_n+i)^{m+1}}X_nh'(0){\rm tr}[\texttt{id}]d\xi_n\sigma(\xi')dx'\nonumber\\
&=\frac{1}{2}Vol(S^{n-2})X_nh'(0)2^m\frac{1}{2}\int_{\Gamma^+}\frac{i\xi_n^2+2m(1-i)\xi_n+2m+i}{(\xi_n-i)^4(\xi_n+i)^{m+1}}d\xi_ndx'\nonumber\\
&=\frac{1}{4}Vol(S^{n-2})X_nh'(0)2^m\frac{2\pi i}{(m+3)!}\left[\frac{i\xi_n^2+2m(1-i)\xi_n+2m+i}{(\xi_n+i)^{m+1}}\right]^{(m+3)}\bigg|_{\xi_n=i}dx'\nonumber\\
&:=Vol(S^{n-2})X_nh'(0)2^m\frac{\pi i}{2(m+3)!}F_0dx'.
\end{align}
where
\begin{align}
F_0&=\left[\frac{i\xi_n^2+2m(1-i)\xi_n+2m+i}{(\xi_n+i)^{m+1}}\right]^{(m+3)}\bigg|_{\xi_n=i}\nonumber\\
&=(2i)^{-2m-4}\bigg(-4iC_{-m-1}^{m+1}+4(m+mi-i)C_{-m-1}^{m+2}+2(2m+mi-i)C_{-m-1}^{m+3}\bigg)(m+3)!.\nonumber\\
\end{align}
\noindent  {\bf (case b-IV)}~$r=-1,~l=-(2m-1),~k=j=|\alpha|=0$.\\
\noindent By (\ref{c2}), we get
\begin{align}\label{c19}
\widetilde{\Psi}_4&=-i\int_{|\xi'|=1}\int^{+\infty}_{-\infty}{\rm tr} [\pi^+_{\xi_n}\sigma_{-1}(\nabla_X^{S(TM)}D^{-1})\times
\partial_{\xi_n}\sigma_{-(2m-1)}(D^{-(2m-1)})](x_0)d\xi_n\sigma(\xi')dx'.
\end{align}
 By (3.24) in \cite{Wa7}, we have
\begin{align}\label{c20}
\partial_{\xi_n}\sigma_{-(2m-1)}(D^{-(2m-1)})(x_0)\bigg|_{|\xi'|=1}
&=\sqrt{-1}\bigg(\frac{c(dx_n)}{(1+\xi_n^2)^{m}}-\frac{2m[\xi_nc(\xi')+\xi_n^2c(dx_n)]}{(1+\xi_n^2)^{m+1}} \bigg).
\end{align}
By Lemma \ref{lem2} and Lemma \ref{lem3}, we have
\begin{align}\label{ggg}
\sigma_{-1}(\nabla_X^{S(TM)}D^{-1})&=\sigma_{1}(\nabla_X^{S(TM)})\sigma_{-2}(D^{-1})+\sigma_{0}(\nabla_X^{S(TM)})\sigma_{-1}(D^{-1})+\sum_{j=1}^n\partial_{\xi_j}\sigma_{1}(\nabla_X^{S(TM)})D_{x_j}[\sigma_{-1}(D^{-1})]\nonumber\\
&:=A^1+A^2+A^3,
\end{align}
where
\begin{align}
A^1(x_0)&=\sqrt{-1}\sum_{j=1}^nX_j\xi_j\bigg[\frac{c(\xi)\sigma_{0}(D)c(\xi)}{|\xi|^4}+\frac{c(\xi)}{|\xi|^6}\sum_jc(dx_j)
\Big(\partial_{x_j}(c(\xi))|\xi|^2-c(\xi)\partial_{x_j}(|\xi|^2)\Big)\bigg];\nonumber\\
A^2(x_0)&=A(X)\frac{\sqrt{-1}c(\xi)}{|\xi|^2};\nonumber\\
A^3(x_0)&=X_n\bigg(\frac{\sqrt{-1}\partial_{x_n}c(\xi')}{|\xi|^2}-\frac{\sqrt{-1}c(\xi)|\xi'|^2h'(0)}{|\xi|^4}\bigg)(x_0).
\end{align}
We note that
\begin{align}\label{45}
\sigma_{0}(D)(x_0)&=-\frac{1}{4}\sum_{s,t,i}\omega_{s,t}(\widetilde{e}_i)
(x_{0})c(\widetilde{e}_i)c(\widetilde{e}_s)c(\widetilde{e}_t):=Q^0(x_0).
\end{align}
Firstly, by $Q^0=c_0c(dx_n)=-\frac{3}{4}h'(0)c(dx_n)$, the following results are obtained by further calculation of $A^1(x_0$)
\begin{align}
A^1(x_0)&=\sqrt{-1}\sum_{j=1}^{n-1}X_j\xi_j\bigg(\frac{3\xi_n^4+4\xi_n^2-7}{4(1+\xi_n^2)^3}h'(0)c(dx_n)+\frac{3\xi_n^2+7\xi_n}{2(1+\xi_n^2)^3}h'(0)c(\xi')+\frac{1}{(1+\xi_n^2)^2}c(\xi')c(dx_n)\partial_{x_n}(c(\xi'))\nonumber\\
&-\frac{\xi_n}{(1+\xi_n^2)^2}\partial_{x_n}(c(\xi'))\bigg)+\sqrt{-1}X_n\bigg(\frac{3\xi_n^5+4\xi_n^3-7\xi_n}{4(1+\xi_n^2)^3}h'(0)c(dx_n)+\frac{3\xi_n^3+7\xi_n^2}{2(1+\xi_n^2)^3}h'(0)c(\xi')\nonumber\\
&+\frac{\xi_n}{(1+\xi_n^2)^2}c(\xi')c(dx_n)\partial_{x_n}(c(\xi'))-\frac{\xi_n^2}{(1+\xi_n^2)^2}\partial_{x_n}(c(\xi'))\bigg).
\end{align}
If we omit some items that have no contribution for computing ${\bf \widetilde{\Psi}_4}$, by the Cauchy integral formula, we obtain
\begin{align}\label{64}
\pi^+_{\xi_n}A^1(x_0)&=\frac{3i\xi_n^2+4\xi_n+i}{8(\xi_n-i)^3}X_nh'(0)c(dx_n)-\frac{-2\xi_n^2+3i\xi_n}{4(\xi_n-i)^3}X_nh'(0)c(\xi')\nonumber\\
&+\frac{1}{4(\xi_n-i)^2}X_nc(\xi')c(dx_n)\partial_{x_n}(c(\xi'))+\frac{i\xi_n}{4(\xi_n-i)^2}X_n\partial_{x_n}(c(\xi')).
\end{align}
By the relation of the Clifford action, we have
\begin{align}\label{ggg32}
&{\rm tr} [\pi^+_{\xi_n}A^1\times
\partial_{\xi_n}\sigma_{-(2m-1)}(D^{-(2m-1)})](x_0)=\frac{3\xi_n^2-5i\xi_n}{8(\xi_n-i)^{m+3}(\xi_n+i)^{m}} X_nh'(0){\rm tr}[\texttt{id}]\nonumber\\
&-\frac{3m\xi_n^4-4mi\xi_n^3+m\xi_n^2}{4(\xi_n-i)^{m+3}(\xi_n+i)^{m}} X_nh'(0){\rm tr}[\texttt{id}]+\frac{2mi\xi_n^4+3m\xi_n^2}{2(\xi_n-i)^{m+4}(\xi_n+i)^{m+1}} X_nh'(0){\rm tr}[\texttt{id}].
\end{align}
Then
 \begin{align}\label{65}
&-i\int_{|\xi'|=1}\int^{+\infty}_{-\infty}{\rm tr} [\pi^+_{\xi_n}A^1\times
\partial_{\xi_n}\sigma_{-(2m-1)}(D^{-(2m-1)})](x_0)d\xi_n\sigma(\xi')dx'\nonumber\\
&=-i\int_{|\xi'|=1}\int^{+\infty}_{-\infty}\frac{3\xi_n^2-5i\xi_n}{8(\xi_n-i)^{m+3}(\xi_n+i)^{m}} X_nh'(0){\rm tr}[\texttt{id}]d\xi_n\sigma(\xi')dx'\nonumber\\
&-i\int_{|\xi'|=1}\int^{+\infty}_{-\infty}-\frac{3m\xi_n^4-4mi\xi_n^3+m\xi_n^2}{4(\xi_n-i)^{m+3}(\xi_n+i)^{m}} X_nh'(0){\rm tr}[\texttt{id}]d\xi_n\sigma(\xi')dx'\nonumber\\
&-i\int_{|\xi'|=1}\int^{+\infty}_{-\infty}\frac{2mi\xi_n^4+3m\xi_n^2}{2(\xi_n-i)^{m+4}(\xi_n+i)^{m+1}} X_nh'(0){\rm tr}[\texttt{id}]d\xi_n\sigma(\xi')dx'\nonumber\\
&=-\frac{i}{8}Vol(S^{n-2})X_nh'(0)2^m\int_{\Gamma^{+}}\frac{3\xi_n^2-5i\xi_n}{(\xi_n-i)^{m+3}(\xi_n+i)^{m}}d\xi_ndx'\nonumber\\
&+\frac{i}{4}Vol(S^{n-2})X_nh'(0)2^m\int_{\Gamma^{+}}\frac{3m\xi_n^4-4mi\xi_n^3+m\xi_n^2}{4(\xi_n-i)^{m+3}(\xi_n+i)^{m}}d\xi_ndx'\nonumber\\
&-\frac{i}{2}Vol(S^{n-2})X_nh'(0)2^m\int_{\Gamma^{+}}\frac{2mi\xi_n^4+3m\xi_n^2}{(\xi_n-i)^{m+4}(\xi_n+i)^{m+1}}d\xi_ndx'\nonumber\\
&=-\frac{i}{8}Vol(S^{n-2})X_nh'(0)2^m\frac{2\pi i}{(m+2)!}\left[\frac{3\xi_n^2-5i\xi_n}{(\xi_n+i)^{m}}\right]^{(m+2)}\bigg|_{\xi_n=i}dx'\nonumber\\
&+\frac{i}{4}Vol(S^{n-2})X_nh'(0)2^m\frac{2\pi i}{(m+2)!}\left[\frac{3m\xi_n^4-4mi\xi_n^3+m\xi_n^2}{(\xi_n+i)^{m}}\right]^{(m+2)}\bigg|_{\xi_n=i}dx'\nonumber\\
&-\frac{i}{2}Vol(S^{n-2})X_nh'(0)2^m\frac{2\pi i}{(m+3)!}\left[\frac{2mi\xi_n^4+3m\xi_n^2}{(\xi_n+i)^{m+1}}\right]^{(m+3)}\bigg|_{\xi_n=i}dx'\nonumber\\
&:=\bigg(\frac{\pi}{4(m+2)!}G_0
-\frac{\pi}{2(m+2)!}G_1
+\frac{\pi }{(m+3)!}G_2\bigg)Vol(S^{n-2})X_nh'(0)2^mdx'\nonumber\\
\end{align}
where
\begin{align}
G_0&=\left[\frac{3\xi_n^2-5i\xi_n}{(\xi_n+i)^{m}}\right]^{(m+2)}\bigg|_{\xi_n=i}=(2i)^{-2m-2}\bigg(-12C_{-m}^m-22C_{-m}^{m+1}-8C_{-m}^{m+2}\bigg)(m+2)!\nonumber\\
G_1&=\left[\frac{3m\xi_n^4-4mi\xi_n^3+m\xi_n^2}{(\xi_n+i)^{m}}\right]^{(m+2)}\bigg|_{\xi_n=i}=-i^{-2m-4}2^{-2m-1}m\bigg(24C_{-m}^{m-2}+(48+16i)C_{-m}^{m-1}+(34+24i)C_{-m}^{m}\nonumber\\
&+(10+12i)C_{-m}^{m+1}+(1+2i)C_{-m}^{m+2}\bigg)(m+2)!;\nonumber\\
G_2&=\left[\frac{2mi\xi_n^4+3m\xi_n^2}{(\xi_n+i)^{m+1}}\right]^{(m+3)}\bigg|_{\xi_n=i}=-i^{-2m-6}2^{-2m-4}\bigg(32mC_{-m-1}^{m-1}+64mC_{-m-1}^{m}+(36m-12)C_{-m-1}^{m+1}\nonumber\\
&+(4m+12)C_{-m-1}^{m+2}+(3-m)C_{-m-1}^{m+3}\bigg)(m+3)!.\nonumber\\
\end{align}
Secondly, for $A^2$, further calculation leads to new results
\begin{align}\label{66}
\pi^+_{\xi_n} A^2(x_0)&=\pi^+_{\xi_n}\Big(A(X)\frac{\sqrt{-1}c(\xi)}{|\xi|^2}\Big)\nonumber\\
&=\frac{1}{2(\xi_n-i)}A(X)c(\xi')+\frac{i}{2(\xi_n-i)}A(X)c(dx_n).
\end{align}
Next
\begin{align}
&{\rm tr} [\pi^+_{\xi_n}A^2\times
\partial_{\xi_n}\sigma_{-(2m-1)}(D^{-(2m-1)})](x_0)\nonumber\\
&=\frac{(-2mi+i)\xi_n^2-2m\xi_n+i}{2(\xi_n-i)^{m+2}(\xi_n+i)^{m+1}}{\rm tr}[A(X)c(\xi')c(dx_n)]+\frac{(2m-1)\xi_n^2-1}{2(\xi_n-i)^{m+2}(\xi_n+i)^{m+1}}{\rm tr}[A(X)c(dx_n)c(dx_n)]\nonumber\\
&-\frac{mi\xi_n}{(\xi_n-i)^{m+2}(\xi_n+i)^{m+1}}{\rm tr}[A(X)c(\xi')c(\xi')].
\end{align}
Because ${\rm tr}[A(X)c(dx_n)c(dx_n)]={\rm tr}[A(X)c(\xi')c(\xi')]=0$, and ${\rm tr }[A(X)c(\xi')c(dx_n)]$ has no contribution for computing $\widetilde{\Psi}_4$.
Then we obtain
\begin{align}
-i\int_{|\xi'|=1}\int^{+\infty}_{-\infty}{\rm tr} [\pi^+_{\xi_n}A^2\times
\partial_{\xi_n}\sigma_{-(2m-1)}(\widetilde{D}^{-(2m-1)})](x_0)d\xi_n\sigma(\xi')dx'=0.
\end{align}
Thirdly, for $A^3$, we get
\begin{align}\label{66pp}
\pi^+_{\xi_n} A^3(x_0)&=\pi^+_{\xi_n}\bigg[X_n\bigg(\frac{\sqrt{-1}\partial_{x_n}c(\xi')}{|\xi|^2}-\frac{\sqrt{-1}c(\xi)|\xi'|^2h'(0)}{|\xi|^4}\bigg)\bigg]\nonumber\\
&=\frac{1}{2(\xi_n-i)}X_n\partial_{x_n}[c(\xi')]+\frac{\xi_n-2i}{4(\xi_n-i)}h'(0)X_nc(\xi')-\frac{1}{4(\xi_n-i)}X_nh'(0)c(dx_n).
\end{align}
Moreover
\begin{align}
{\rm tr} [\pi^+_{\xi_n}A^3\times
\partial_{\xi_n}\sigma_{-(2m-1)}(D^{-(2m-1)})](x_0)=\frac{i\xi_n^2+2m\xi_n+i}{4(\xi_n-i)^{m+3}(\xi_n+i)^{m+1}}X_nh'(0){\rm tr}[\texttt{id}].
\end{align}
Then, we have
\begin{align}
&-i\int_{|\xi'|=1}\int^{+\infty}_{-\infty}{\rm tr} [\pi^+_{\xi_n}A^3\times
\partial_{\xi_n}\sigma_{-(2m-1)}(\widetilde{D}^{-(2m-1)})](x_0)d\xi_n\sigma(\xi')dx'\nonumber\\
&-i\int_{|\xi'|=1}\int^{+\infty}_{-\infty}\frac{i\xi_n^2+2m\xi_n+i}{4(\xi_n-i)^{m+3}(\xi_n+i)^{m+1}}X_nh'(0){\rm tr}[\texttt{id}]d\xi_n\sigma(\xi')dx'\nonumber\\
&=-iV0l(S^{n-2})X_nh'(0)2^m\int_{\Gamma^+}\frac{i\xi_n^2+2m\xi_n+i}{4(\xi_n-i)^{m+3}(\xi_n+i)^{m+1}}d\xi_ndx'\nonumber\\
&=-\frac{i}{4}V0l(S^{n-2})X_nh'(0)2^m\frac{2\pi i}{(m+2)!}\left[\frac{i\xi_n^2+2m\xi_n+i}{(\xi_n+i)^{m+1}}\right]^{(m+2)}\bigg|_{\xi_n=i}dx'\nonumber\\
&:=V0l(S^{n-2})X_nh'(0)2^m\frac{\pi}{2(m+2)!}G_3dx',
\end{align}
where
\begin{align}
G_3=\left[\frac{i\xi_n^2+2m\xi_n+i}{(\xi_n+i)^{m+1}}\right]^{(m+2)}\bigg|_{\xi_n=i}=i^{-2m-4}4^{-m-1}\bigg(2C_{-m-1}^m+(2m+10)C_{-m-1}^{m+1}-mC_{-m-1}^{m+2}\bigg)(m+2)!.
\end{align}
Therefore, we obtain
\begin{align}\label{6666}
\widetilde{\Psi}_4&=-i\int_{|\xi'|=1}\int^{+\infty}_{-\infty}{\rm tr} [\pi^+_{\xi_n}(A^1+A^2+A^3)\times
\partial_{\xi_n}\sigma_{-(2m-1)}(D^{-(2m-1)})](x_0)d\xi_n\sigma(\xi')dx'\nonumber\\
&=\bigg(-\frac{i}{8}\frac{2\pi i}{(m+2)!}G_0
+\frac{i}{4}\frac{2\pi i}{(m+2)!}G_1
-\frac{i}{2}\frac{2\pi i}{(m+3)!}G_2-\frac{i}{4}\frac{2\pi i}{(m+2)!}G_3\bigg)V0l(S^{n-2})X_nh'(0)2^mdx'.
\end{align}
\noindent {\bf (case b-V)}~$r=0,~\ell=-2m,~k=j=|\alpha|=0$.\\
By (\ref{c2}), we get
\begin{align}\label{c23}
\widetilde{\Psi}_5&=-i\int_{|\xi'|=1}\int^{+\infty}_{-\infty}{\rm tr} [\pi^+_{\xi_n}\sigma_{0}(\nabla_X^{S(TM)}D^{-1})\times
\partial_{\xi_n}\sigma_{-2m}(D^{-(2m-1)})](x_0)d\xi_n\sigma(\xi')dx'\nonumber\\
&=i\int_{|\xi'|=1}\int^{+\infty}_{-\infty}{\rm tr} [\partial_{\xi_n}\pi^+_{\xi_n}\sigma_{0}(\nabla_X^{S(TM)}D^{-1})\times
\sigma_{-2m}(D^{-(2m-1)})](x_0)d\xi_n\sigma(\xi')dx'.
\end{align}
By (3.37) in \cite{Wa7}, we have
\begin{align}\label{c24}
\sigma_{-2m}(D^{-(2m-1)})(x_0)
&=\frac{(-2m-1)h'(0)c(\mathrm{d}x_n)}{4(1+\xi_n^2)^{\frac{n}{2}}}-2m\xi_n(1+\xi_n^2)^{-m-1}\partial_{x_n}(c(\xi')(x_0)+mi(1+\xi_n^2)^{-m+1}\nonumber\\
&[c(\xi')+\xi_nc(dx_n)]\times\bigg[\frac{-ih'(0)c(\xi')c(dx_n)-(2m+1)ih'(0)c(\xi')}{2(1+\xi_n^2)^2}-\frac{2ih'(0)\xi_n}{(1+\xi_n^2)^3}\bigg]\nonumber\\
&-[c(\xi')+\xi_nc(dx_n)]h'(0)\xi_n[m^2+m][(1+\xi_n^2)^{-m-2}].\nonumber\\
\end{align}
By Lemma \ref{lem2} and Lemma \ref{lem3}, we have
\begin{align}\label{c25}
\partial_{\xi_n}\pi^+_{\xi_n}\sigma_{0}(\nabla_X^{S(TM)}D^{-1})
&=\frac{i}{2(\xi_n-i)^2}\sum_{j=1}^{n-1}X_j\xi_j
c(\xi')+\frac{1}{2(\xi_n-i)^2}\sum_{j=1}^{n-1}X_j\xi_j
c(dx_n)\nonumber\\
&+\frac{1}{2(\xi_n-i)^2}X_nc(\xi')+\frac{i}{2(\xi_n-i)^2}X_nc(dx_n).
\end{align}
When we omit some items that have no contribution for computing ${\bf \widetilde{\Psi}_5}$. Then, we have
\begin{align}\label{cc55}
&{\rm tr} [\partial_{\xi_n}\pi^+_{\xi_n}\sigma_{0}(\nabla_X^{S(TM)}D^{-1})\times
\sigma_{-2m}(D^{-(2m-1)})](x_0)\nonumber\\
&=\frac{(2m^2-m)i\xi_n^2}{8(\xi_n-i)^{m+3}(\xi_n+i)^{m+1}}X_nh'(0){\rm tr}[\texttt{id}]+\frac{(im^2+m^2+m-3mi)\xi_n}{2(\xi_n-i)^{m+4}(\xi_n+i)^{m+2}}X_nh'(0){\rm tr}[\texttt{id}]\nonumber\\
&+\frac{(2m^2+3m+1)i\xi_n^2}{8(\xi_n-i)^{m+3}(\xi_n+i)^{m+1}}X_nh'(0){\rm tr}[\texttt{id}]+\frac{(3m+1)i}{8(\xi_n-i)^{m+3}(\xi_n+i)^{m+1}}X_nh'(0){\rm tr}[\texttt{id}].
\end{align}
Therefore, we have
\begin{align}\label{c27}
\widetilde{\Psi}_5&=i\int_{|\xi'|=1}\int^{+\infty}_{-\infty}{\rm tr} [\partial_{\xi_n}\pi^+_{\xi_n}\sigma_{0}(\nabla_X^{S(TM)}D^{-1})\times
\sigma_{-2m}(D^{-(2m-1)})](x_0)d\xi_n\sigma(\xi')dx'\nonumber\\
&=i\int_{|\xi'|=1}\int^{+\infty}_{-\infty}\frac{(2m^2-m)i\xi_n^2}{8(\xi_n-i)^{m+3}(\xi_n+i)^{m+1}}X_nh'(0){\rm tr}[\texttt{id}]d\xi_n\sigma(\xi')dx'\nonumber\\
&+i\int_{|\xi'|=1}\int^{+\infty}_{-\infty}\frac{(im^2+m^2+m-3mi)\xi_n}{2(\xi_n-i)^{m+4}(\xi_n+i)^{m+2}}X_nh'(0){\rm tr}[\texttt{id}]d\xi_n\sigma(\xi')dx'\nonumber\\
&+i\int_{|\xi'|=1}\int^{+\infty}_{-\infty}\frac{(2m^2+3m+1)i\xi_n^2}{8(\xi_n-i)^{m+3}(\xi_n+i)^{m+1}}X_nh'(0){\rm tr}[\texttt{id}]d\xi_n\sigma(\xi')dx'\nonumber\\
&+i\int_{|\xi'|=1}\int^{+\infty}_{-\infty}\frac{(3m+1)i}{8(\xi_n-i)^{m+3}(\xi_n+i)^{m+1}}X_nh'(0){\rm tr}[\texttt{id}]d\xi_n\sigma(\xi')dx'\nonumber\\
&=\frac{(2m^2-m)i}{8}V0l(S^{n-2})X_nh'(0)2^m\frac{2\pi i}{(m+2)!}\left[\frac{\xi_n}{(\xi_n+i)^{m+1}}\right]^{(m+2)}\bigg|_{\xi_n=i}dx'\nonumber\\
&+\frac{(i-1)m^2+mi+3m}{2}V0l(S^{n-2})X_nh'(0)2^m\frac{2\pi i}{(m+3)!}\left[\frac{\xi_n}{(\xi_n+i)^{m+2}}\right]^{(m+3)}\bigg|_{\xi_n=i}dx'\nonumber\\
&-\frac{2m^2+3m+1}{8}V0l(S^{n-2})X_nh'(0)2^m\frac{2\pi i}{(m+2)!}\left[\frac{\xi_n^2}{(\xi_n+i)^{m+1}}\right]^{(m+2)}\bigg|_{\xi_n=i}dx'\nonumber\\
&-\frac{3m+1}{2}V0l(S^{n-2})X_nh'(0)2^m\frac{2\pi i}{(m+2)!}\left[\frac{1}{(\xi_n+i)^{m+1}}\right]^{(m+2)}\bigg|_{\xi_n=i}dx'\nonumber\\
&:=\bigg(-\frac{(2m^2-m)\pi}{4(m+2)!}H_0
+\frac{((i-1)m^2+mi+3m)\pi i}{(m+3)!}H_1
-\frac{(2m^2+3m+1)\pi i}{4(m+2)!}H_2\nonumber\\
&-\frac{(3m+1)\pi i}{(m+2)!}H_3\bigg)V0l(S^{n-2})X_nh'(0)2^mdx'.
\end{align}
where
\begin{align}
H_0&=\left[\frac{\xi_n}{(\xi_n+i)^{m+1}}\right]^{(m+2)\bigg|_{\xi_n=i}}=-i^{-2m-4}2^{-2m-3}\bigg(2C_{-m-1}^{m+1}+C_{-m-1}^{m+2}\bigg)(m+2)!;\nonumber\\
H_1&=\left[\frac{\xi_n}{(\xi_n+i)^{m+2}}\right]^{(m+3)\bigg|_{\xi_n=i}}=(2i)^{-2m-4}\bigg(2iC_{-m-1}^{m+2}+iC_{-m-1}^{m+3}\bigg)(m+3)!;\nonumber\\
H_2&=\left[\frac{\xi_n^2}{(\xi_n+i)^{m+1}}\right]^{(m+2)\bigg|_{\xi_n=i}}=-(2i)^{-2m-3}\bigg(4C_{-m-1}^{m}+4C_{-m-1}^{m+1}+C_{-m-1}^{m+2}\bigg)(m+2)!;\nonumber\\
H_3&=\left[\frac{1}{(\xi_n+i)^{m+1}}\right]^{(m+2)\bigg|_{\xi_n=i}}=(2i)^{-2m-3}A_{-m-1}^{m+2}.\nonumber\\
\end{align}
Now $\widetilde{\Psi}$ is the sum of the {\bf  (case b-I)}-{\bf (case b-V)}. Therefore, we get
\begin{align}\label{c28}
\widetilde{\Psi}&=\sum_{i=1}^5\widetilde{\Psi}_i\nonumber\\
&=\bigg\{-\frac{1}{4}Vol(S^{n-2})2^m\partial_{x_n}(X_n)\frac{2\pi i}{(m+2)!}E_0
+\bigg(-\frac{1}{4}\frac{2\pi i}{(m+3)!}E_1+\frac{1}{4}\frac{2\pi i}{(m+3)!}F_0-\frac{i}{8}\frac{2\pi i}{(m+3)!}G_0\nonumber\\
&-\frac{i}{4}\frac{2\pi i}{(m+2)!}G_1+\frac{(2m^2-m)i}{8}\frac{2\pi i}{(m+2)!}H_0
+\frac{(i-1)m^2+mi+3m}{2}\frac{2\pi i}{(m+3)!}H_1\nonumber\\
&
-\frac{2m^2+3m+1}{8}\frac{2\pi i}{(m+2)!}H_2
-\frac{3m+1}{2}\frac{2\pi i}{(m+2)!}H_3\bigg)V0l(S^{n-2})X_nh'(0)2^m\bigg\}dx'.
\end{align}
Then, by (\ref{a2}) and (\ref{c28}), we obtain following theorem
\begin{thm}\label{cthmb1}
Let $M$ be an $n=2m+1$-dimensional oriented
compact spin manifold with boundary $\partial M$, then we get the following equality:
\begin{align}
\label{cb2773}
&\widetilde{{\rm Wres}}[\pi^+(\nabla_X^{S(TM)}D^{-1})\circ\pi^+(D^{-(2m-1)})]\nonumber\\
&=\int_{\partial M}\bigg\{-\frac{1}{4}Vol(S^{n-2})2^m\partial_{x_n}(X_n)\frac{2\pi i}{(m+2)!}E_0
+\bigg(-\frac{\pi i}{2(m+3)!}E_1+\frac{\pi i}{2(m+3)!}F_0+\frac{\pi }{4(m+3)!}G_0\nonumber\\
&+\frac{\pi }{2(m+2)!}G_1-\frac{(2m^2-m)\pi}{4(m+2)!}H_0
+\frac{((i-1)m^2+mi+3m)\pi i}{(m+3)!}H_1-\frac{(2m^2+3m+1)\pi i}{4(m+2)!}H_2\nonumber\\
&-\frac{(3m+1)\pi i}{(m+2)!}H_3\bigg)Vol(S^{n-2})X_nh'(0)2^m\bigg\}d{\rm Vol_{M}}
\end{align}
\end{thm}

\section{The noncommutative residue $\widetilde{{\rm Wres}}[\pi^+(\nabla_X^{S(TM)}D^{-2})\circ\pi^+(D^{-(2m-2)})]$ on odd dimensional manifolds with boundary}
\label{section:5}
In this section, we compute the noncommutative residue $\widetilde{{\rm Wres}}[\pi^+(\nabla_X^{S(TM)}D^{-2})\circ\pi^+(D^{-(2m-2)})]$ on $2m+1$ dimensional oriented
compact spin manifolds with boundary.

Similar to \cite{Wa3}, by Proposition \ref{prop1}, we can get
\begin{align}
\label{b14}
&\widetilde{{\rm Wres}}[\pi^+(\nabla_X^{S(TM)}D^{-2})\circ\pi^+(D^{-(2m-2)})]=\int_{\partial M}\widehat{\Psi},
\end{align}
where
\begin{align}
\label{b1115}
\widehat{\Psi}&=\int_{|\xi'|=1}\int^{+\infty}_{-\infty}\sum^{\infty}_{j, k=0}\sum\frac{(-i)^{|\alpha|+j+k+1}}{\alpha!(j+k+1)!}
\times {\rm tr}_{\wedge^*T^*M\bigotimes\mathbb{C}}[\partial^j_{x_n}\partial^\alpha_{\xi'}\partial^k_{\xi_n}\sigma^+_{r}(\nabla_X^{S(TM)}D^{-2})(x',0,\xi',\xi_n)
\nonumber\\
&\times\partial^\alpha_{x'}\partial^{j+1}_{\xi_n}\partial^k_{x_n}\sigma_{l}(D^{-(2m-2)})(x',0,\xi',\xi_n)]d\xi_n\sigma(\xi')dx',
\end{align}
and the sum is taken over $r+l-k-j-|\alpha|-1=-(2m+1),~~r\leq -1,~~l\leq -(2m-2)$.

Next, we need to compute $\int_{\partial M} \widehat{\Psi}$. The sum is taken over $
r+l-k-j-|\alpha|=-2m,~~r\leq -1,~~l\leq -(2m-2),$ then we have the following five cases:

\noindent  {\bf (case c-I)}~$r=-1,~l=-(2m-2),~k=j=0,~|\alpha|=1$.\\
\noindent By (\ref{b1115}), we get
\begin{equation}
\label{b24}
\widehat{\Psi}_1=-\int_{|\xi'|=1}\int^{+\infty}_{-\infty}\sum_{|\alpha|=1}
 {\rm tr}[\partial^\alpha_{\xi'}\pi^+_{\xi_n}\sigma_{-1}(\nabla_X^{S(TM)}D^{-2})\times
 \partial^\alpha_{x'}\partial_{\xi_n}\sigma_{-(2m-2)}(D^{-(2m-2)})](x_0)d\xi_n\sigma(\xi')dx'.
\end{equation}
By Lemma 2.2 in \cite{Wa3}, for $i<n$, then
\begin{equation}
\label{b25}
\partial_{x_i}\sigma_{-(2m-2)}(D^{-(2m-2)})(x_0)=\partial_{x_i}{(|\xi|^{(2-2m)})}(x_0)
=\partial_{x_i}(|\xi|^2)^{(1-m)}(x_0)=(1-m) (|\xi|^2)^{-m}
 \partial_{x_i}(|\xi|^2)(x_0)=0,
\end{equation}
\noindent so $\widehat{\Psi}_1=0$.\\
 \noindent  {\bf (case c-II)}~$r=-1,~l=-(2m-2),~k=|\alpha|=0,~j=1$.\\
\noindent By (\ref{b1115}), we get
\begin{equation}
\label{b26}
\widehat{\Psi}_2=-\frac{1}{2}\int_{|\xi'|=1}\int^{+\infty}_{-\infty} {\rm
tr} [\partial_{x_n}\pi^+_{\xi_n}\sigma_{-1}(\nabla_X^{S(TM)}D^{-2})\times
\partial_{\xi_n}^2\sigma_{-(2m-2)}(D^{-(2m-2)})](x_0)d\xi_n\sigma(\xi')dx'.
\end{equation}
\noindent By (3.16) in \cite{Wa6}, we have
\begin{align}\label{b237}
\partial^2_{\xi_n} \sigma_{-(2m-2)}(D^{-(2m-2)})(x_0)
&=\partial^2_{\xi_n}\big((|\xi|^2)^{1-m}\big) (x_0)
=\partial_{\xi_n}\Big( (1-m) (|\xi|^{2})^{-\frac{n}{2}}  \partial_{\xi_n} (|\xi|^2)\Big) (x_0)\nonumber\\
&=(1-m)(-m)(|\xi|^{2})^{-m-1} \big(\partial_{\xi_n}|\xi|^2\big)^{2}(x_0)+(1-m)(|\xi|^{2})^{-m} \partial^2_{\xi_n}(|\xi|^2 (x_0))\nonumber\\
&=\Big((4m-2) \xi_n^{2}-2 \Big)(m-1)(1+\xi_n^{2})^{(-m-1)}.
\end{align}
By Lemma \ref{lem3}, we have
\begin{align}\label{b27}
\partial_{x_n}\sigma_{-1}(\nabla_X^{S(TM)}D^{-2})&=\partial_{x_n}\bigg(i\sum_{j=1}^nX_j\xi_j|\xi|^{-2}\bigg)\nonumber\\
&=i\sum_{j=1}^{n-1}\xi_j\frac{\partial X_j}{\partial x_n}\frac{1}{1+\xi_n^2}+i\frac{\partial X_n}{\partial x_n}\frac{\xi_n}{1+\xi_n^2}-i\sum_{j=1}^{n-1}\frac{X_j\xi_jh'(0)|\xi'|^2}{(1+\xi_n^2)}-iX_nh'(0)|\xi'|^2\frac{\xi_n}{(1+\xi_n^2)}.
\end{align}
Moreover
\begin{align}\label{b28}
&\partial_{x_n}\pi^+_{\xi_n}\sigma_{-1}(\nabla_X^{S(TM)}D^{-2})\nonumber\\
&=\pi^+_{\xi_n}\partial_{x_n}\sigma_{-1}(\nabla_X^{S(TM)}D^{-2})\nonumber\\
&=\frac{1}{2(\xi_n-i)}\sum_{j=1}^{n-1}\xi_j\frac{\partial X_j}{\partial x_n}+\frac{i}{2(\xi_n-i)}\frac{\partial X_n}{\partial x_n}+\frac{2i-\xi_n}{4(\xi_n-i)^2}\sum_{j=1}^{n-1}X_j\xi_jh'(0)|\xi'|^2-\frac{1}{4(\xi_n-i)^2}X_nh'(0)|\xi'|^2.
\end{align}
Then, we have
\begin{align}\label{35}
&{\rm
tr} [\partial_{x_n}\pi^+_{\xi_n}\sigma_{-1}(\nabla_X^{S(TM)}D^{-2})\times
\partial_{\xi_n}^2\sigma_{-(2m-2)}(D^{-(2m-2)})](x_0)\nonumber\\
&=(m-1)\sum_{j=1}^{n-1}\xi_j\frac{\partial X_j}{\partial x_n}\frac{(2m-1)\xi_n^2-1}{(\xi_n+i)^{m+1}(\xi_n-i)^{m+2}}{\rm tr}[\texttt{id}]+(m-1)\frac{\partial X_n}{\partial x_n}\frac{(2m-1)i\xi_n^2-i}{(\xi_n+i)^{m+1}(\xi_n-i)^{m+2}}{\rm tr}[\texttt{id}]\nonumber\\
&+(m-1)\sum_{j=1}^{n-1}X_j\xi_jh'(0)|\xi'|^2\frac{((2m-1)\xi_n^2-1)(2i-\xi_n)}{(\xi_n+i)^{m+1}(\xi_n-i)^{m+3}}{\rm tr}[\texttt{id}]+(m-1)X_nh'(0)|\xi'|^2\frac{(2m-1)\xi_n^2-1}{(\xi_n+i)^{m+1}(\xi_n-i)^{m+3}}{\rm tr}[\texttt{id}].
\end{align}
We omit some items that have no contribution for computing ${\bf \widehat{\Psi}_2}$. Then, we get
\begin{align}\label{35kkk}
\widehat{\Psi}_2&=-\frac{1}{2}\int_{|\xi'|=1}\int^{+\infty}_{-\infty} {\rm
tr} [\partial_{x_n}\pi^+_{\xi_n}\sigma_{-1}(\nabla_X^{S(TM)}D^{-2})\times
\partial_{\xi_n}^2\sigma_{-(2m-2)}(D^{-(2m-2)})](x_0)d\xi_n\sigma(\xi')dx'\nonumber\\
&=-\frac{1}{2}\int_{|\xi'|=1}\int^{+\infty}_{-\infty}(m-1)\frac{\partial X_n}{\partial x_n}\frac{(2m-1)i\xi_n^2-i}{(\xi_n+i)^{m+1}(\xi_n-i)^{m+2}}{\rm tr}[\texttt{id}]d\xi_n\sigma(\xi')dx'\nonumber\\
&-\frac{1}{2}\int_{|\xi'|=1}\int^{+\infty}_{-\infty}(m-1)X_nh'(0)|\xi'|^2\frac{(2m-1)\xi_n^2-1}{(\xi_n+i)^{m+1}(\xi_n-i)^{m+3}}{\rm tr}[\texttt{id}]d\xi_n\sigma(\xi')dx'\nonumber\\
&=-\frac{1}{2}(m-1)Vol(S^{n-2})2^m\frac{\partial X_n}{\partial x_n}\int_{\Gamma^+}\frac{(2m-1)i\xi_n^2-i}{(\xi_n+i)^{m+1}(\xi_n-i)^{m+2}}d\xi_ndx'\nonumber\\
&-\frac{1}{2}(m-1)Vol(S^{n-2})X_nh'(0)2^m|\xi'|^2\int_{\Gamma^+}\frac{(2m-1)\xi_n^2-1}{(\xi_n+i)^{m+1}(\xi_n-i)^{m+3}}d\xi_ndx'\nonumber\\
&=-\frac{1}{2}(m-1)Vol(S^{n-2})2^m\frac{\partial X_n}{\partial x_n}\frac{2\pi i}{(m+1)!}\left[\frac{(2m-1)i\xi_n^2-i}{(\xi_n+i)^{m+1}}\right]^{(m+1)}\bigg|_{\xi_n=i}dx'\nonumber\\
&+\frac{1}{4}(m-1)Vol(S^{n-2})2^mX_nh'(0)|\xi'|^2\frac{2\pi i}{(m+2)!}\left[\frac{(2m-1)\xi_n^2-1}{(\xi_n+i)^{m+1}}\right]^{(m+2)}\bigg|_{\xi_n=i}dx'\nonumber\\
&:=(m-1)\bigg(-\frac{\partial X_n}{\partial x_n}\frac{\pi i}{(m+1)!}I_0+X_nh'(0)\frac{\pi i}{2(m+2)!}I_1\bigg)Vol(S^{n-2})2^mdx',
\end{align}
where
\begin{align}
I_0&=\left[\frac{(2m-1)i\xi_n^2-i}{(\xi_n+i)^{m+1}}\right]^{(m+1)}\bigg|_{\xi_n=i}=-(2i)^{-2m-1}\bigg((4m-2)C_{-m-1}^{m-1}+(4m-2)C_{-m-1}^m+mC_{-m-1}^{m+1}\bigg)(m+1)!\nonumber\\
I_1&=\left[\frac{(2m-1)\xi_n^2-1}{(\xi_n+i)^{m+1}}\right]^{(m+2)}\bigg|_{\xi_n=i}=(2i)^{-2m-3}\bigg((4-8m)C_{-m-1}^{m}+(4-8m)C_{-m-1}^{m+1}-2mC_{-m-1}^{m+2}\bigg)(m+2)!.
\end{align}
\noindent  {\bf (case c-III)}~$r=-1,~l=-(2m-2),~j=|\alpha|=0,~k=1$.\\
\noindent By (\ref{b1115}), we get
\begin{align}\label{36}
\widehat{\Psi}_3&=-\frac{1}{2}\int_{|\xi'|=1}\int^{+\infty}_{-\infty}
{\rm tr} [\partial_{\xi_n}\pi^+_{\xi_n}\sigma_{-1}(\nabla_X^{S(TM)}D^{-2})\times
\partial_{\xi_n}\partial_{x_n}\sigma_{-(2m-2)}(D^{-(2m-2)})](x_0)d\xi_n\sigma(\xi')dx'\nonumber\\
&=\frac{1}{2}\int_{|\xi'|=1}\int^{+\infty}_{-\infty}
{\rm tr} [\partial_{\xi_n}^2\pi^+_{\xi_n}\sigma_{-1}(\nabla_X^{S(TM)}D^{-2})\times
\partial_{x_n}\sigma_{-(2m-2)}(D^{-(2m-2)})](x_0)d\xi_n\sigma(\xi')dx'.
\end{align}
By (3.21) in \cite{Wa6}, we have
\begin{align}\label{37}
\partial_{x_n} \big(\sigma_{-(2m-2)}(D^{-(2m-2)})\big)(x_0)
=\partial_{x_n}\big((|\xi|^2)^{1-m}\big) (x_0)
=h'(0)(1-m)(1+\xi_{n}^{2})^{-m}.
\end{align}
By Lemma \ref{lem2} and Lemma \ref{lem3}, we have
\begin{align}\label{aa38}
\pi^+_{\xi_n}\sigma_{-1}(\nabla_X^{S(TM)}D^{-2})=\frac{1}{2(\xi_n-i)}\sum_{j=1}^{n-1}X_j\xi_j+\frac{i}{2(\xi_n-i)}X_n.
\end{align}
By further calculation, we have
\begin{align}\label{mmmmm}
\partial_{\xi_n}^2\pi^+_{\xi_n}\sigma_{-1}(\nabla_X^{S(TM)}D^{-2})=\frac{1}{(\xi_n-i)^3}\sum_{j=1}^{n-1}X_j\xi_j+\frac{i}{(\xi_n-i)^3}X_n.
\end{align}
Moreover
\begin{align}
&{\rm tr} [\partial_{\xi_n}^2\pi^+_{\xi_n}\sigma_{-1}(\nabla_X^{S(TM)}D^{-2})\times
\partial_{x_n}\sigma_{-(2m-2)}(D^{-(2m-2)})](x_0)\nonumber\\
&=(1-m)h'(0)\sum_{j=1}^{n-1}X_j\xi_j\frac{1}{(\xi_n-i)^{m+3}(\xi_n+i)^{m}}{\rm tr}[\texttt{id}]+(1-m)X_nh'(0)\frac{i}{(\xi_n-i)^{m+3}(\xi_n+i)^{m}}{\rm tr}[\texttt{id}].
\end{align}
Next, we perform the corresponding integral calculation on the above results. Then, we get
\begin{align}\label{41}
\widehat{\Psi}_3&=\frac{1}{2}\int_{|\xi'|=1}\int^{+\infty}_{-\infty}
{\rm tr} [\partial_{\xi_n}^2\pi^+_{\xi_n}\sigma_{-1}(\nabla_X^{S(TM)}D^{-2})\times
\partial_{x_n}\sigma_{-(2m-2)}(D^{-(2m-2)})](x_0)d\xi_n\sigma(\xi')dx'\nonumber\\
&=\frac{1}{2}\int_{|\xi'|=1}\int^{+\infty}_{-\infty}
(1-m)X_nh'(0)\frac{i}{(\xi_n-i)^{m+3}(\xi_n+i)^m}{\rm tr}[\texttt{id}]d\xi_n\sigma(\xi')dx'\nonumber\\
&=\frac{1}{2}(1-m)Vol(S^{n-2})X_nh'(0)2^m\int_{\Gamma^+}\frac{i}{(\xi_n-i)^{m+3}(\xi_n+i)^m}d\xi_ndx'\nonumber\\
&=\frac{1}{2}(1-m)Vol(S^{n-2})X_nh'(0)2^m\frac{2\pi i}{(m+2)!}\left[\frac{i}{(\xi_n+i)^m}\right]^{(m+2)}\bigg|_{\xi_n=i}dx'\nonumber\\
&:=(1-m)Vol(S^{n-2})X_nh'(0)2^m\frac{\pi i}{2(m+2)!}J_0dx'\nonumber\\
\end{align}
where
\begin{align}
J_0&=\left[\frac{i}{(\xi_n+i)^m}\right]^{(m+2)}\bigg|_{\xi_n=i}=-i^{-2m-3}2^{-2m-1}A^{m+2}_m.\nonumber\\
\end{align}
\noindent  {\bf (case c-IV)}~$r=-1,~l=-(2m-1),~k=j=|\alpha|=0$.\\
\noindent By (\ref{b1115}), we get
\begin{align}\label{42}
\widehat{\Psi}_4&=-i\int_{|\xi'|=1}\int^{+\infty}_{-\infty}{\rm tr} [\pi^+_{\xi_n}\sigma_{-1}(\nabla_X^{S(TM)}D^{-2})\times
\partial_{\xi_n}\sigma_{-(2m-1)}(D^{-(2m-2)})](x_0)d\xi_n\sigma(\xi')dx'\nonumber\\
&=i\int_{|\xi'|=1}\int^{+\infty}_{-\infty}{\rm tr} [\partial_{\xi_n}\pi^+_{\xi_n}\sigma_{-1}(\nabla_X^{S(TM)}D^{-2})\times
\sigma_{-(2m-1)}(D^{-(2m-2)})](x_0)d\xi_n\sigma(\xi')dx'.
\end{align}
By (3.30) in \cite{Wa6}, we have
\begin{align}\label{43}
\sigma_{-(2m-1)}(D^{-(2m-2)})
&=-\frac{2m^2-m-1}{2}h'(0)\frac{i\xi_n}{(1+\xi_n^2)^m}-(m-1)h'(0)\frac{2i\xi_n}{(1+\xi_n^2)^{m+1}}\nonumber\\
&-(m^2-3m+2)h'(0)\frac{i\xi_n}{(1+\xi_n^2)^{m+1}}.
\end{align}
 By Lemma \ref{lem2} and Lemma \ref{lem3}, we have
\begin{align}\label{45}
\partial_{\xi_n}\pi^+_{\xi_n}\sigma_{-1}(\nabla_X^{S(TM)}D^{-2})&=-\frac{1}{2(\xi_n-i)^2}\sum_{j=1}^{n-1}X_j\xi_j-\frac{i}{2(\xi_n-i)^2}X_.
\end{align}
Moreover
\begin{align}
&{\rm tr} [\partial_{\xi_n}\pi^+_{\xi_n}\sigma_{-1}(\nabla_X^{S(TM)}D^{-2})\times
\sigma_{-(2m-1)}(D^{-(2m-2)})](x_0)\nonumber\\
&=\frac{(2m^2-m-1)i\xi_n}{4(\xi_n+i)^m(\xi_n-i)^{m+2}}h'(0)\sum_{j=1}^{n-1}X_j\xi_j{\rm tr}[\texttt{id}]-\frac{(2m^2-m-1)\xi_n}{4(\xi_n+i)^m(\xi_n-i)^{m+2}}X_nh'(0){\rm tr}[\texttt{id}]\nonumber\\
&+\frac{(m-1)i\xi_n}{(\xi_n+i)^m(\xi_n-i)^{m+2}})h'(0)\sum_{j=1}^{n-1}X_j\xi_j{\rm tr}[\texttt{id}]-\frac{(m-1)\xi_n}{(\xi_n+i)^m(\xi_n-i)^{m+2}})X_nh'(0){\rm tr}[\texttt{id}]\nonumber\\
&+\frac{(m^2-3m+2)i\xi_n}{2(\xi_n+i)^m(\xi_n-i)^{m+3}}h'(0)\sum_{j=1}^{n-1}X_j\xi_j{\rm tr}[\texttt{id}]-\frac{(m^2-3m+2)\xi_n}{2(\xi_n+i)^m(\xi_n-i)^{m+3}}X_nh'(0){\rm tr}[\texttt{id}].
\end{align}
We omit some items that have no contribution for computing ${\bf \widehat{\Psi}_4}$. Then, we have
\begin{align}\label{39}
\widehat{\Psi}_4&=-i\int_{|\xi'|=1}\int^{+\infty}_{-\infty}\frac{(2m^2-m-1)\xi_n}{4(\xi_n+i)^m(\xi_n-i)^{m+2}}X_nh'(0){\rm tr}[\texttt{id}]d\xi_n\sigma(\xi')dx'\nonumber\\
&-i\int_{|\xi'|=1}\int^{+\infty}_{-\infty}\frac{(m-1)\xi_n}{(\xi_n+i)^m(\xi_n-i)^{m+2}})X_nh'(0){\rm tr}[\texttt{id}]d\xi_n\sigma(\xi')dx'\nonumber\\
&-i\int_{|\xi'|=1}\int^{+\infty}_{-\infty}\frac{(m^2-3m+2)\xi_n}{2(\xi_n+i)^m(\xi_n-i)^{m+3}}X_nh'(0){\rm tr}[\texttt{id}]d\xi_n\sigma(\xi')dx'\nonumber\\
&=-\frac{(2m^2+3m-5)i}{4}Vol(S^{n-2})X_nh'(0)2^m\int_{\Gamma^+}\frac{\xi_n}{(\xi_n+i)^m(\xi_n-i)^{m+2}}d\xi_ndx'\nonumber\\
&-\frac{(m^2-3m+2)i}{2}Vol(S^{n-2})X_nh'(0)2^m\int_{\Gamma^+}\frac{\xi_n}{(\xi_n+i)^m(\xi_n-i)^{m+3}}d\xi_ndx'\nonumber\\
&=-\frac{(2m^2+3m-5)i}{4}Vol(S^{n-2})X_nh'(0)2^m\frac{2\pi i}{(m+1)!}\left[\frac{\xi_n}{(\xi_n+i)^m}\right]^{(m+1)}\bigg|_{\xi_n=i}dx'\nonumber\\
&-\frac{(m^2-3m+2)i}{2}Vol(S^{n-2})X_nh'(0)2^m\frac{2\pi i}{(m+2)!}\left[\frac{\xi_n}{(\xi_n+i)^m}\right]^{(m+2)}\bigg|_{\xi_n=i}dx'\nonumber\\
&:=\bigg(\frac{(2m^2+3m-5)\pi }{2(m+1)!}K_0+\frac{(m^2-3m+2)\pi}{(m+2)!}K_1\bigg)Vol(S^{n-2})X_nh'(0)2^mdx',
\end{align}
where
\begin{align}
K_0&=\left[\frac{\xi_n}{(\xi_n+i)^m}\right]^{(m+1)}\bigg|_{\xi_n=i}=(2i)^{-2m-1}\bigg(2iC_{-m}^{m}+iC_{-m}^{m+1}\bigg)(m+1)!\nonumber\\
 &=\nonumber\\
 K_1&=\left[\frac{\xi_n}{(\xi_n+i)^m}\right]^{(m+2)}\bigg|_{\xi_n=i}=(2i)^{-2m-2}\bigg(2iC_{-m}^{m+1}+iC_{-m}^{m+2}\bigg)(m+2)!.
\end{align}
\noindent {\bf  (case c-V)}~$r=-2,~\ell=-(2m-2),~k=j=|\alpha|=0$.\\
By (\ref{b1115}), we get
\begin{align}\label{61}
\widehat{\Psi}_5=-i\int_{|\xi'|=1}\int^{+\infty}_{-\infty}{\rm tr} [\pi^+_{\xi_n}\sigma_{-2}(\nabla_X^{S(TM)}D^{-2})\times
\partial_{\xi_n}\sigma_{-(2m-2)}(D^{-(2m-2)})](x_0)d\xi_n\sigma(\xi')dx'.
\end{align}
By (3.33) in \cite{Wa6}, we have
\begin{equation}\label{62}
\partial_{\xi_n}\sigma_{-(2m-2)}(D^{-(2m-2)})(x_0)
=\partial_{\xi_n}((|\xi|^2)^{1-m})(x_0)=2(1-m)\xi_n(1+\xi_n^2)^{-m}.
\end{equation}
By Lemma \ref{lem2} and Lemma \ref{lem3}, we have
\begin{align}\label{673}
\sigma_{-2}(\nabla^{S(TM)}_XD^{-2})(x_0)&=\sigma_{0}(\nabla^{S(TM)})\sigma_{-2}(D^{-2})+\sigma_{1}(\nabla^{S(TM)})\sigma_{-3}(D^{-2})+\sum_{j=1}^n\partial_{\xi_j}\sigma_{1}(\nabla^{S(TM)})D_{x_j}[\sigma_{-2}(D^{-2})]\nonumber\\
&:=B^1+B^2+B^3,
\end{align}
where
\begin{align}
B^1(x_0)&=A(X)|\xi|^{-2};\nonumber\\
B^2(x_0)&=\sqrt{-1}\sum_{j=1}^nX_j\xi_j\bigg[\frac{i}{2(1+\xi_n^2)^2}
   h'(0)\sum_{k<n}\xi_kc(\widetilde{e}_k)c(\widetilde{e}_n)-\frac{5i\xi_n^3+9i\xi_n}{2(1+\xi_n^2)^3}h'(0)\bigg];\nonumber\\
B^3(x_0)&=-X_n\frac{h'(0)|\xi'|^2}{|\xi|^4}.
\end{align}
Firstly, the following results are obtained by further calculation of $B^1(x_0)$
\begin{align}
\pi^+_{\xi_n} B^1(x_0)=\frac{i}{2(\xi_n-i)}A(X).
\end{align}
Then
\begin{align}
{\rm tr} [\pi^+_{\xi_n} B^1\times
\partial_{\xi_n}\sigma_{-(2m-2)}(D^{-(2m-2)})](x_0)=-\frac{(1-m)i\xi_n}{(\xi_n-i)(1+\xi_n^2)^2}{\rm tr}[A(X)].
\end{align}
We note that ${\rm tr}[A(X)]=0$, then
 \begin{align}\label{65}
-i\int_{|\xi'|=1}\int^{+\infty}_{-\infty}{\rm tr} [\pi^+_{\xi_n}B^1\times
\partial_{\xi_n}\sigma_{-(2m-2)}(\widetilde{D}^{-(2m-2)})](x_0)d\xi_n\sigma(\xi')dx'=
0.
\end{align}
Secondly, for $B^2(x_0)$, further calculation leads to new results
\begin{align}\label{66}
\pi^+_{\xi_n} B^2(x_0)&=\frac{2+i\xi_n}{8(\xi_n-i)^2}h'(0)\sum_{j=1}^{n-1}X_j\xi_j\sum_{k<n}\xi_kc(\widetilde{e}_k)c(\widetilde{e}_n)+\frac{i}{8(\xi_n-i)^2}X_nh'(0)\sum_{k<n}\xi_kc(\widetilde{e}_k)c(\widetilde{e}_n)\nonumber\\
&-\frac{3i\xi_n+4}{4(\xi_n-i)^3}h'(0)\sum_{j=1}^{n-1}X_j\xi_j-\frac{3i\xi_n^2+4\xi_n}{4(\xi_n-i)^3}X_nh'(0).
\end{align}
By $\sum_{j=k<n}^{n-1}{\rm tr} [c(\widetilde{e}_k)c(\widetilde{e}_n)]=0$, and we omit some items that have no contribution for computing ${\bf \widehat{\Psi}_5}$.
\begin{align}
{\rm tr} [\pi^+_{\xi_n}B^2\times
\partial_{\xi_n}\sigma_{-(2m-2)}(\widetilde{D}^{-(2m-2)})](x_0)=-(1-m)\frac{3i\xi_n^3+4\xi_n^2}{2(\xi_n-i)^{m+3}(\xi_n+i)^m}X_nh'(0){\rm tr}[\texttt{id}].
\end{align}
Thirdly, for $B^3(x_0)$, we get
\begin{align}\label{661pp}
\pi^+_{\xi_n} B^3(x_0)
&=\frac{2+i\xi_n}{4(\xi_n-i)^2}X_nh'(0).
\end{align}
Then
\begin{align}
&{\rm tr} [\pi^+_{\xi_n}B^3\times
\partial_{\xi_n}\sigma_{-(2m-2)}(D^{-(2m-2)})](x_0)=(1-m)\frac{i\xi_n^2+2\xi_n}{2(\xi_n-i)^2(1+\xi_n^2)^m}X_nh'(0){\rm tr}[\texttt{id}].
\end{align}
Moreover
\begin{align}
&{\rm tr} [\pi^+_{\xi_n}(B^2+B^3)\times
\partial_{\xi_n}\sigma_{-(2m-2)}(D^{-(2m-2)})](x_0)=\frac{(1-m)(2i\xi_n^3+\xi_n^2+2i\xi_n)}{2(\xi_n-i)^3(1+\xi_n^2)^m}X_nh'(0){\rm tr}[\texttt{id}].
\end{align}
Then, we have
\begin{align}
\widehat{\Psi}_5&=-i\int_{|\xi'|=1}\int^{+\infty}_{-\infty}{\rm tr} [\pi^+_{\xi_n}(B^1+B^2+B^3\times
\partial_{\xi_n}\sigma_{-(2m-2)}(D^{-(2m-2)})](x_0)d\xi_n\sigma(\xi')dx'\nonumber\\
&=-i\int_{|\xi'|=1}\int^{+\infty}_{-\infty}\frac{(1-m)(2i\xi_n^3+\xi_n^2+2i\xi_n)}{2(\xi_n-i)^3(1+\xi_n^2)^m}X_nh'(0){\rm tr}[\texttt{id}]d\xi_n\sigma(\xi')dx'\nonumber\\
&=\frac{(1-m)i}{2}Vol(S^{n-2})X_nh'(0)2^m\int_{\Gamma^+} \frac{2i\xi_n^3+\xi_n^2+2i\xi_n}{(\xi_n-i)^{m+3}(\xi_n+i)^m}d\xi_ndx'\nonumber\\
&=\frac{(1-m)i}{2}Vol(S^{n-2})X_nh'(0)2^m\frac{2\pi i}{(m+2)!}\left[\frac{2i\xi_n^3+\xi_n^2+2i\xi_n}{(\xi_n+i)^m}\right]^{(m+2)}\bigg|_{\xi_n=i}dx'\nonumber\\
&=(m-1)Vol(S^{n-2})X_nh'(0)2^m\frac{\pi}{(m+2)!}L_0dx',
\end{align}
where
\begin{align}
L_0=\left[\frac{2i\xi_n^3+\xi_n^2+2i\xi_n}{(\xi_n+i)^m}\right]^{(m+2)}\bigg|_{\xi_n=i}=(2i)^{-2m-2}\bigg(-16C_{-m}^{m+1}+20C_{-m}^{m}+4C_{-m}^{m+1}-C_{-m}^{m+2}\bigg)(m+2)!.
\end{align}
Now $\widehat{\Psi}$ is the sum of the {\bf  (case c-I)}-{\bf  (case c-V)}. Therefore, we get
\begin{align}\label{795}
\widehat{\Psi}&=\sum_{i=1}^5\widehat{\Psi}_i\nonumber\\
&=\bigg\{-(m-1)Vol(S^{n-2})2^m\frac{\partial X_n}{\partial x_n}\frac{\pi i}{(m+1)!}I_0+\bigg((m-1)\frac{\pi i}{2(m+2)!}I_1
+(1-m)\frac{\pi i}{(m+2)!}J_0\nonumber\\
&+\frac{(2m^2+3m-5)\pi }{2(m+1)!}K_0+\frac{(m^2-3m+2)\pi }{(m+2)!}K_1+(m-1)\frac{2\pi i}{(m+2)!}L_0\bigg)Vol(S^{n-2})X_nh'(0)2^m\bigg\}dx' .
\end{align}
Then, by (\ref{795}), we obtain following theorem
\begin{thm}\label{thmb1}
Let $M$ be an $n=2m+1$-dimensional oriented
compact spin manifold with boundary $\partial M$, then we get the following equality:
\begin{align}
\label{b2773}
&\widetilde{{\rm Wres}}[\pi^+(\nabla_X^{S(TM)}D^{-2})\circ\pi^+(D^{-(2m-2)})]\nonumber\\
&=\int_{\partial M} \bigg\{-(m-1)Vol(S^{n-2})2^m\frac{\partial X_n}{\partial x_n}\frac{\pi i}{(m+1)!}I_0+\bigg((m-1)\frac{\pi i}{2(m+2)!}I_1
+(1-m)\frac{\pi i}{(m+2)!}J_0\nonumber\\
&+\frac{(2m^2+3m-5)\pi }{2(m+1)!}K_0+\frac{(m^2-3m+2)\pi }{(m+2)!}K_1+(m-1)\frac{2\pi i}{(m+2)!}L_0\bigg)Vol(S^{n-2})X_nh'(0)2^m\bigg\}d{\rm Vol_{M}}.
\end{align}
\end{thm}
\section*{Acknowledgements}
This work was supported by NSFC No.11771070. And the partial research of the third author was
supported by DUFE202159 and Basic research Project of the Education Department of Liaoning Province (Grant No. LJKQZ20222442). The authors thank the referee for his (or her) careful reading and helpful comments.

\section*{}


\begin{thebibliography}{00}
\bibitem{Ac} Ackermann T. A note on the Wodzicki residue. J. Geom.Phys. 1996, 20: 404-406.
\bibitem{Co1} Connes A. Quantized calculus and applications. 11th International Congress of Mathematical Physics(Paris,1994), Internat Press, Cambridge, MA, 1995, 15-36.
\bibitem{Co2} Connes A. The action functinal in Noncommutative geometry. Comm. Math. Phys. 1998, 117: 673-683.
\bibitem{DL} Dabrowski L, Sitarz A, Zalecki P. Spectral Metric and Einstein Functionals. Adv. Math. 2023, 427: 109128.
\bibitem{FGLS} Fedosov B V, Golse F, Leichtnam E, Schrohe E. The noncommutative residue for manifolds with boundary. J. Funct. Anal. 1996, 142: 1-31.
\bibitem{Gu} Guillemin V W. A new proof of Weyl's formula on the asymptotic distribution of eigenvalues. Adv. Math. 1985, 55(2): 131-160.
\bibitem{KW} Kalau W, Walze M. Gravity, Noncommutative geometry and the Wodzicki residue. J. Geom. Physics. 1995, 16: 327-344.
\bibitem{Ka} Kastler D. The Dirac Operator and Gravitation. Comm. Math. Phys. 1995, 166: 633-643.
\bibitem{Wa5} Wang J, Wang Y. The Kastler-Kalau-Walze type theorem for six-dimensional manifolds with boundary. J. Math. Phys. 2015, 56: 052501.
\bibitem{Wa6} Wang J, Wang Y. A general Kastler-Kalau-Walze type theorem for manifolds with boundary. Int. J. Geom. Methods M. Phys. 2016,
13 (1): 1650003
\bibitem{Wa1} Wang Y. Diffential forms and the Wodzicki residue for Manifolds with Boundary. J. Geom. Physics. 2006, 56: 731-753.
\bibitem{Wa2} Wang Y. Diffential forms the Noncommutative Residue for Manifolds with Boundary in the non-product Case. Lett. Math. Phys. 2006, 77: 41-51.
\bibitem{Wa3} Wang Y. Gravity and the Noncommutative Residue for Manifolds with Boundary. Lett. Math. Phys. 2007, 80: 37-56.
\bibitem{Wa4} Wang Y. Lower-Dimensional Volumes and Kastler-kalau-Walze Type Theorem for Manifolds with Boundary. Commun. Theor. Phys. 2010, 54: 38-42.
\bibitem{Wa7} Wang Y. General Kastler-Kalau-Walze type theoremsfor manifolds with boundary II. Int. J. Geom. Methods M. Phys. 2019, 16(2): 1950028.
\bibitem{Wo} Wodzicki M. local invariants of spectral asymmetry. Invent. Math. 1995, 75(1): 143-178.
 \bibitem{Wu1} Wu T, Wang Y. A general Dabrowski-Sitarz-Zalecki type theorems for manifold with boundary. ArXiv:2308.15850.
\bibitem{Wu2} Wu T, Wang Y. The generalized noncommutative residue and the Kastler-Kalau-Walze type theorem . ArXiv.2204.11021.
\bibitem{Y} Yu Y. The Index Theorem and The Heat Equation Method, Nankai Tracts in Mathematics-Vol.2, World Scientific Publishing, 2001.
\end{thebibliography}
\end{document}